\documentclass[12pt]{article}

\usepackage[hyperref]{package_RFb}

\usepackage[top=2cm,bottom=2cm,left=3cm,right=3cm]{geometry}
\usepackage[numbers]{natbib}
\usepackage{authblk}
\usepackage{stmaryrd}
\usepackage{caption}

\title{Demographic inference for spatially heterogeneous populations using long shared haplotypes}
\author[1]{Raphaël Forien\footnote{raphael.forien@inrae.fr}}
\author[2]{Harald Ringbauer\footnote{harald{\textunderscore}ringbauer@eva.mpg.de}}
\author[3]{Graham Coop\footnote{gmcoop@ucdavis.edu}}

\affil[1]{INRAE - BioSP, Centre INRAE PACA, 228~route de l'aérodrome, Domaine St-Paul - Site Agroparc, 84914, Avignon Cedex 9, France}
\affil[2]{Department of Archaeogenetics, Max Planck Institute for Evolutionary Anthropology, Deutscher Platz~6, 04103, Leipzig, Germany}
\affil[3]{Center for Population Biology, Department of Evolution and Ecology, University of California, 2320~Storer Hall, CA~95616, Davis, United States}

\graphicspath{{figures/}}
\usepackage[export]{adjustbox}
\usepackage{titlesec}

\NewDocumentCommand{\PM}{}{\lbrace +, - \rbrace}
\NewDocumentCommand{\Block}{}{B}

\NewDocumentCommand\Etheta{o m g}{
	\mathbb{E}^\theta%
	\IfNoValueF{#1}{_{#1}}
	\left[%
	\IfNoValueF{#3}{\left.}%
	#2%
	\IfNoValueF{#3}{\:\right|\: #3}
	\right]}

\begin{document}

\maketitle

    \begin{abstract}
        We introduce a modified spatial $\Lambda$-Fleming-Viot process to model the ancestry of individuals in a population occupying a continuous spatial habitat divided into two areas by a sharp discontinuity of the dispersal rate and effective population density. We derive an analytical formula for the expected number of shared haplotype segments between two individuals depending on their sampling locations.
        This formula involves the transition density of a skew diffusion which appears as a scaling limit of the ancestral lineages of individuals in this model. We then show that this formula can be used to infer the dispersal parameters and the effective population density of both regions, using a composite likelihood approach, and we demonstrate the efficiency of this method on a range of simulated data sets.

        \textbf{Keywords:} population genetics, spatial $\Lambda$-Fleming-Viot process, spatial coalescent, segments of shared haplotypes, skew Brownian motion, isolation by distance.
    \end{abstract}

\section{Introduction}

Within spatially structured populations, genetic variation is shaped by a range of evolutionary processes, such as dispersal, demography, mutations, recombination, and selection.
When the average distance travelled by individuals over their lifetime is much shorter than the whole range of the population, spatial patterns of neutral genetic diversity can arise, in which genetic similarity between individuals decreases with the geographic distance separating them \cite{sharbel_genetic_2000, novembre_genes_2008, aguillon_deconstructing_2017, battey_space_2020}. 
This pattern, called ``isolation by distance'', is also predicted by theoretical models of reproducing populations occupying a spatially structured habitat and undergoing limited dispersal \cite{wright_isolation_1943,barton_neutral_2002}.

Early studies focused on how geographic distance affects the probability that a given pair of sampled individuals inherit the same allele from a common ancestor at some given locus (called the probability of identity by descent) \cite{malecot_les_1948,kimura_stepping_1964,sawyer_results_1976}, assuming that the population density and migration patterns of individuals are homogeneous across the whole region occupied by the population.
These theoretical predictions can then be compared to empirically observed pairwise genetic distance statistics in order to infer demographic parameters relating to dispersal and population density \cite{rousset_genetic_1997,barton_inference_2013}.
The two main parameters governing the dynamics of these models are the effective dispersal rate, which corresponds to the standard deviation of the displacement from an individual's birthplace to that of their parent, and Wright's neighbourhood size, which can be seen as an effective density of reproducing individuals over one generation (it is inversely proportional to the rate of coalescence of nearby lineages).

The assumptions of the above models, however, are often violated on large spatial and temporal scales, as many populations present significant large-scale heterogeneity. In those situations, efficient inference methods suited to spatially heterogeneous populations are required.

One possible approach is to rely on computationally heavy simulations of complex demographic scenarios using one of the available population models (see for example \cite{guindon_demographic_2016}). Another promising method, MAPS, has been introduced in \cite{al-asadi_estimating_2019}.
Its theoretical foundation is close to that of the present paper, but MAPS infers a large number of parameters, corresponding to the dispersal rates and neighbourhood size at each point of a spatial grid, and needs to implement a strongly regularising penalisation to avoid overfitting.

Here, we take a different approach in which we consider that large-scale variations in dispersal and population density can be summed up by dividing space into a small number of connected regions in which the population behaves according to the homogeneous model, with different parameters in each subregion.
More precisely, here we study a model in which the effective dispersal rate and population density change abruptly when crossing a narrow interface, and derive theoretical predictions regarding the genetic similarity of two individuals depending on their sampling locations (whether on the same side of the interface or not). We then adapt previous inference methods to estimate the parameters of the model from a set of genetic samples sampled on both sides of the interface.

To achieve this, we make use of the distribution of long continuous segments of (near) identical haplotypes between pairs of individuals, so-called identity by descent segments (or IBD segments). Over the past decade, the use of IBD segments has emerged as a fruitful tool to study recent ancestry in human populations \cite{palamara_length_2012,carmi_renewal_2014, baharian_great_2016,fournier_haplotype-based_2022}.
These segments are shared between pairs of individuals, having been inherited, unbroken by recombination, from a recent common ancestor.  Since recombination breaks up the genome at each generation, thus decoupling the genealogies at different loci, long (\textit{i.e.} > 2-4 cM) continuous segments of genome shared identical by descent between individuals are the result of recent coancestry (typically less than 60 generations ago \cite{ralph_geography_2013}). Hence the geographic structure of these long segments carries information about the recent demography of the sample. The expected length of these IBD segments is inversely proportional to the age of the common ancestor from which they are inherited, thus providing information on the timing of coalescence events in the genealogy of the sample.

These IBD segments need to be inferred computationally from sets of genetic sequences (either whole genome sequences or dense single nucleotide data). Several methods have been developed to detect long IBD segments \cite{browning_fast_2011,browning_probabilistic_2020,zhou_fast_2020}, and it is also possible to quantify the uncertainty of the resulting estimates (see \cite{ralph_geography_2013,browning_probabilistic_2020}).

For example, Ralph and Coop \cite{ralph_geography_2013} analysed the European POPRES dataset \cite{nelson_population_2008} to detect IBD segments in a sample of 2,257 individuals from around Europe. They identified a total of 1.9 million IBD segments, finding that most Europeans, even separated by several thousand kilometres, share genetic and genealogical ancestors who lived less than 1,000 years ago. Importantly, they found that the number of shared IBD segments between individuals decreases approximately exponentially with the geographic distance separating them.

Patterns of identity by descent sharing across geographic space allow the inference of dispersal distances in the recent past \cite{al-asadi_estimating_2019}.
Such an inference method was developed in \cite{ringbauer_inferring_2017}, and then applied to the IBD segment calls for the Central and Eastern European subregion \cite{ralph_geography_2013}.
\cite{ringbauer_inferring_2017} focused on Central and Eastern Europe because the whole European sample does not present a homogeneous decrease in the number of IBD segments with geographic distance. A plausible explanation for this spatial heterogeneity is that Western and Eastern Europe broadly present different dispersal rates or effective neighbourhood sizes \cite{ralph_geography_2013}.

Here, we present a new method to infer demographic parameters when both the dispersal and neighbourhood sizes are discontinuous across a straight interface.
Our method is based on a theoretical analysis of a model constructed in the framework of the spatial $\Lambda$-Fleming-Viot process introduced in \cite{barton_new_2010}.
In this model, the population occupies a continuous spatial habitat, and reproduction takes place at a random sequence of events, during which a fraction of the population in a bounded region is replaced by the offspring of some randomly chosen (pair of) individual(s).

We model the heterogeneous dispersal and neighbourhood size by letting the radius and the fraction of replaced individuals depend on the region in which the events take place.
This model can be seen as a continuous space generalisation of Nagylaki's stepping stone model with variable migration \cite{nagylaki_clines_1976}, and indeed we expect the two models to display similar IBD segment sharing patterns over large spatial scales.

We derive the theoretical prediction of the average number of IBD segments of a given length found among two sampled individuals in this model under some assumptions (see Section~\ref{sec:model}), and we use this prediction in the composite likelihood pipeline developed in \cite{ringbauer_inferring_2017} to obtain a software that can jointly infer the dispersal and neighbourhood size parameters in each half-space.
We then apply this method to simulated data to show that it can reliably estimate the parameters used for each simulation in various realistic settings.

The main step in the proof of our analytical result is to show that the rescaled trajectories of ancestral lineages in our model converge in distribution to trajectories of a diffusion with a singular drift at the interface. This type of diffusion is called skew Brownian motion, and was introduced in \cite{ito_brownian_1963} (see also \cite{walsh_diffusion_1978,harrison_skew_1981} and \cite{lejay_constructions_2006} for a review of its properties).
The proof of our result is based on similar arguments as the one established for rescaled random walks in \cite{iksanov_functional_2016}, with several adaptations to deal with a continuous space setting.

The rest of the paper is laid out as follows.
In Section~\ref{sec:model}, we present the spatial $\Lambda$-Fleming-Viot model with heterogeneous dispersal and heterogeneous neighbourhood size, we define the corresponding ancestral process which describes the genome-wide genealogy of a sample of individuals in this model, and we state our main results on the expected number of shared IBD segments.
In Section~\ref{sec:method}, we present the statistical method used for inference, extending the maximum likelihood method of \cite{ringbauer_inferring_2017}.
Section~\ref{sec:results} presents the results of the application of the method to simulated data, and is followed by a discussion in Section~\ref{sec:discussion}.
The proofs of the mathematical results used to compute the expected sharing of IBD segments are provided in the appendix.

\section{Model and theoretical results} \label{sec:model}

\subsection{Model and main assumptions}

We consider a population occupying a continuous two-dimensional spatial habitat, which we identify with $\R^2$. The population evolves through a random sequence of reproduction events, following the framework of the spatial $\Lambda$-Fleming-Viot process (or SLFV for short) introduced in \cite{barton_new_2010} which simulates neutral evolution, forwards in time, of a two-dimensional population. 
At each reproduction event, a fraction of the population in a small disc around the location of the event dies and is replaced by the offspring of a randomly chosen pair of individuals (who are chosen from the same disc).

In this model, the fraction of the local population that is replaced during a reproduction event is a measure of the (inverse) local population size, as it determines the size of the offspring of one pair of individuals relative to the rest of the population. This fraction, also called the impact parameter, is usually denoted by $u$ (see below). The size of the discs affected by reproduction events also measures the strength of dispersal in the population, as it determines the typical distance between parents and their offspring.

In order to account for the spatial heterogeneity of the local population density and of the strength of dispersal, we introduce a modification of the usual spatial $\Lambda$-Fleming-Viot process in which we split $\R^2$ into the positive and negative half spaces, denoted here as $\mathbb{H}_-$ (on the left) and $\mathbb{H}_+$ (on the right). 

We assume that the impact parameter and the radius of each reproduction event depend on the half space in which its centre falls, so that an event whose centre lies in $\mathbb{H}_+$ has impact $u_+$ and radius $r_+$, while an event whose centre falls in $\mathbb{H}_-$ has impact $u_-$ and radius $r_-$.
Note that, in this model, an individual sitting in $\mathbb{H}_+$ (resp. $\mathbb{H}_-$) dies at rate $\pi u_+ (r_+)^2$ (resp. $\pi u_- (r_-)^2$) times the rate of either type of events (individuals who are sitting close to the boundary between the two half-spaces may die at slightly different rates, but we ignore this).
The rates of reproduction events in both half spaces are thus adjusted to ensure that individuals on both sides have the same average lifetime.

\begin{figure}[ht]
    \centering
    \includegraphics[width=0.5\textwidth]{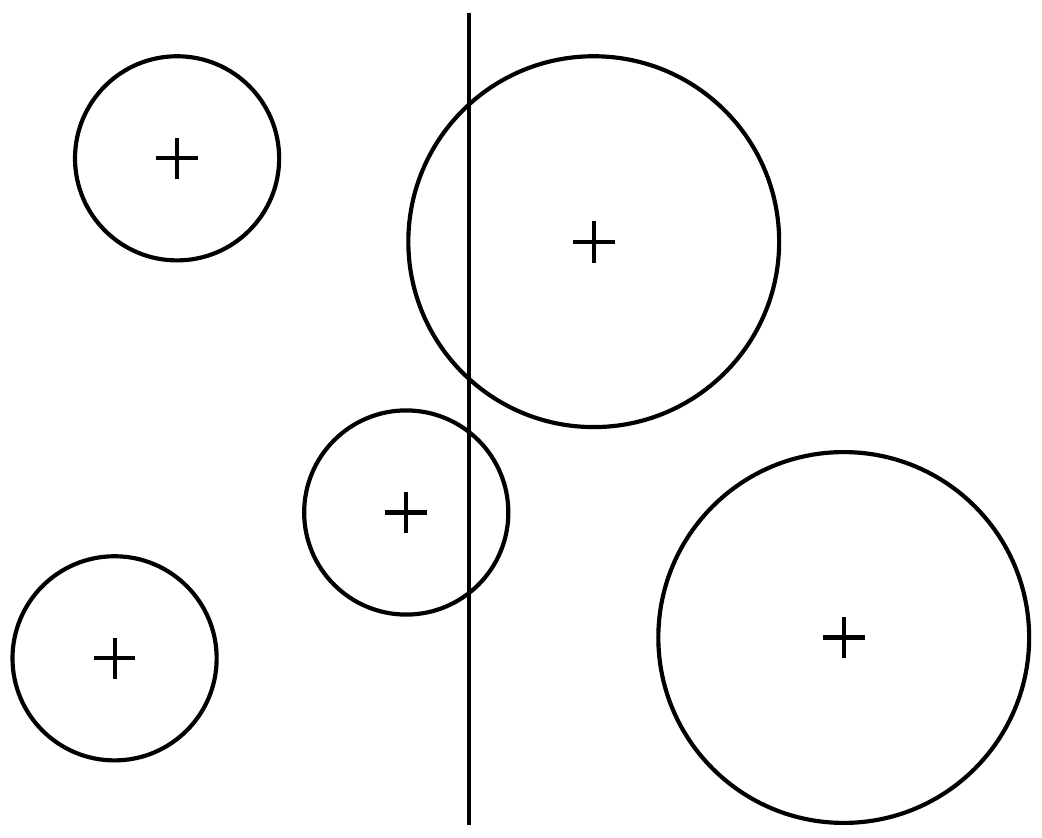}
    \caption{Regions affected by different reproduction events. The vertical line denotes the interface between $\mathbb{H}_-$ (on the left) and $\mathbb{H}_+$ (on the right). As can be seen from the figure, individuals lying in $\mathbb{H}_-$ close enough to $\mathbb{H}_+$ may be involved in reproduction events whose centre falls in $\mathbb{H}_+$ (which allows their offspring to cross from $\mathbb{H}_-$ to $\mathbb{H}_+$ if they are chosen as a parent). The rates of the two types of events are scaled so that individuals that do not live close to the interface are affected at the same rate in both half-spaces.}
    \label{fig:slfv_hetero}
\end{figure}

We focus on a single chromosome and assume that each individual is haploid, i.e. they carry a single copy of this chromosome.
During a reproduction event, each offspring's genome is the product of the random recombination of two randomly selected parental genomes.
On the time scales of interest this should be an accurate approximation, since the genealogies of different chromosomes (homologous or not) can effectively be considered independent because Mendelian segregation ensures that, a few generations in the past, ancestral chromosomes are carried by different individuals with high probability.
We approximate recombination by a Poisson process with rate 1 along the genome, not modelling crossover interference (see Figure~\ref{fig:recombination} of the supplementary material), so that we measure genetic distances along the genome in Morgans. 

\begin{definition} \label{def:recombination_pattern}
A recombination pattern is a random piece-wise constant map $\mathcal{R} : [0,G] \to \lbrace 1, 2 \rbrace$, such that $\mathcal{R}(0)$ equals 1 with probability $1/2$ (and 2 otherwise), and $\mathcal{R}$ changes value at the points of a standard Poisson process on $[0,G]$.
\end{definition}

We then consider a sample of individuals from this population, chosen at random from a fixed set of locations, and we record the IBD segments longer than some threshold length among pairs of individuals in the sample. We give a formal definition of IBD segments in this model below (Definition~\ref{def:IBD_segment}).
We then derive an analytical formula approximating the expected number of such segments for a pair of individuals, as a function depending on their sampling locations (see \eqref{eq:WM_E_N_geq_L}).
This formula is valid provided the sampling locations are sufficiently far apart, the local population density is sufficiently large (or equivalently the impact parameters in each half-space are sufficiently small), and the length of IBD segments is sufficiently short (see the more detailed statement in Appendix~\ref{app:Wright_Malecot}).

\subsection{The ancestral recombination graph of the SLFV with heterogeneous dispersal}
We now formally define a process describing the genealogy of a sample of $n$ individual chromosomes from the model described above. The appropriate tool for this is the so-called ancestral recombination graph, introduced in \cite{hudson_properties_1983} (see also \cite{griffiths_ancestral_1997}), with lineages labelled by their spatial position.
At some time $t$ before the sampling time, the genetic ancestors of the sample can be described by a set of labelled lineages,
\begin{equation*}
    \mathcal{A}_t := \left\lbrace ( \Block^{i}_t, \xi^{i}_t ), 1 \leq i \leq N_t \right\rbrace,
\end{equation*}
where $\xi^{i}_t \in \R^2$ are the locations of the genetic ancestors, and the $\Block^i_t$ are maps from $[0,G]$ to the set of finite subsets of $\lbrace 1, \ldots, n \rbrace$, such that, for any locus $\ell \in [0, G]$, $ \lbrace \Block^{i}_t(\ell), 1 \leq i \leq N_t \rbrace$ is a partition of $\lbrace 1, \ldots, n \rbrace$, possibly with some empty blocks.
Thus, the genetic ancestor, at locus $\ell$ and at time $t$ before sampling, of the $j$-th sampled individual is located at $\xi^i_t$, where $i$ is such that $j \in \Block^i_t(\ell)$.
Note that time runs backwards here, i.e. $\mathcal{A}_t$ describes the genetic ancestors of the sampled individuals $t$ generations in the past.
We assume that $\mathcal{A}_t$ is started from
\begin{equation*}
    \mathcal{A}_0 = \left\lbrace ( \lbrace i \rbrace, x_i ), 1 \leq i \leq n \right\rbrace,
\end{equation*}
where $\lbrace x_i, 1 \leq i \leq n \rbrace$ are the sampling locations.
The genealogy at a specific locus $\ell \in [0,G]$ can then be obtained as
\begin{equation*}
    \mathcal{A}^\ell_t := \left\lbrace ( \Block^i_t(\ell), \xi^i_t), 1 \leq i \leq N_t, \Block^i_t(\ell) \neq \emptyset \right\rbrace.
\end{equation*}

Let us now give the law of $(\mathcal{A}_t, t\geq 0)$ under the model described above.
Fix $u_\alpha \in (0,1]$ and $r_\alpha > 0$ for $\alpha \in \PM$.
For $r > 0$, $x \in \R^2$, the ball of centre $x$ and radius $r$ is denoted by $B(x,r)$ and its volume by $V_r = \pi r^2$.
Also, for $\alpha \in \PM$, define the positive (resp. negative) half space $\mathbb{H}_+$ (resp. $\mathbb{H}_-$) as
\begin{equation*}
    \mathbb{H}_\alpha := \lbrace x \in \R^2 : \alpha x^{(1)} > 0 \rbrace.
\end{equation*}

\begin{definition} \label{def:ARG}
Let $\lbrace x_1, \ldots, x_n \rbrace$ denote the sampling positions, and set, for $1 \leq i \leq n$,
\begin{equation*}
    \xi^{i}_0 = x_i, \quad \Block^{i}_0(\ell) = \lbrace i \rbrace, \quad \forall \ell \in [0,G].
\end{equation*}
For $\alpha \in \PM$, let $\Pi_\alpha$ be a Poisson random measure on $\R_+ \times \mathbb{H}_\alpha$ with intensity measure $(u_\alpha V_{r_\alpha})^{-1} dt \otimes dx$.
For each $(t,x) \in \Pi_\alpha$, a reproduction event with impact $u_\alpha$ takes place in $B(x,r_\alpha)$ following the steps below.
\begin{itemize}
    \item Pick two locations $y_1$ and $y_2$ uniformly from $B(x,r_\alpha)$, independently of each other.
    \item Each lineage located in $B(x,r_\alpha)$ just before the event is involved with probability $u_\alpha$, independently of each other. Let $I$ denote the (random) set of involved lineages.
    \item Each involved lineage $i \in I$ picks a random recombination pattern $\mathcal{R}_i$ following Definition~\ref{def:recombination_pattern}.
    \item All involved lineages are replaced by two lineages $(y_1, \hat{\Block}^1)$ and $(y_2, \hat{\Block}^2)$, where
    \begin{equation*}
        \hat{\Block}^i(\ell) = \bigcup_{j \in I : \mathcal{R}_j(\ell) = i} \Block^j_{t^-}(\ell), \quad i \in \lbrace 1, 2 \rbrace.
    \end{equation*}
\end{itemize}
In addition, if, after a reproduction event, a lineage $i$ is such that $\Block^i_t(\ell) = \emptyset$ for all $\ell \in [0,G]$, we drop this lineage from $\mathcal{A}_t$ (since it corresponds to an ancestor who did not pass down any genetic material to the sampled individuals).
\end{definition}

\begin{figure}[htp]
    \centering
    \includegraphics[width=0.95\textwidth]{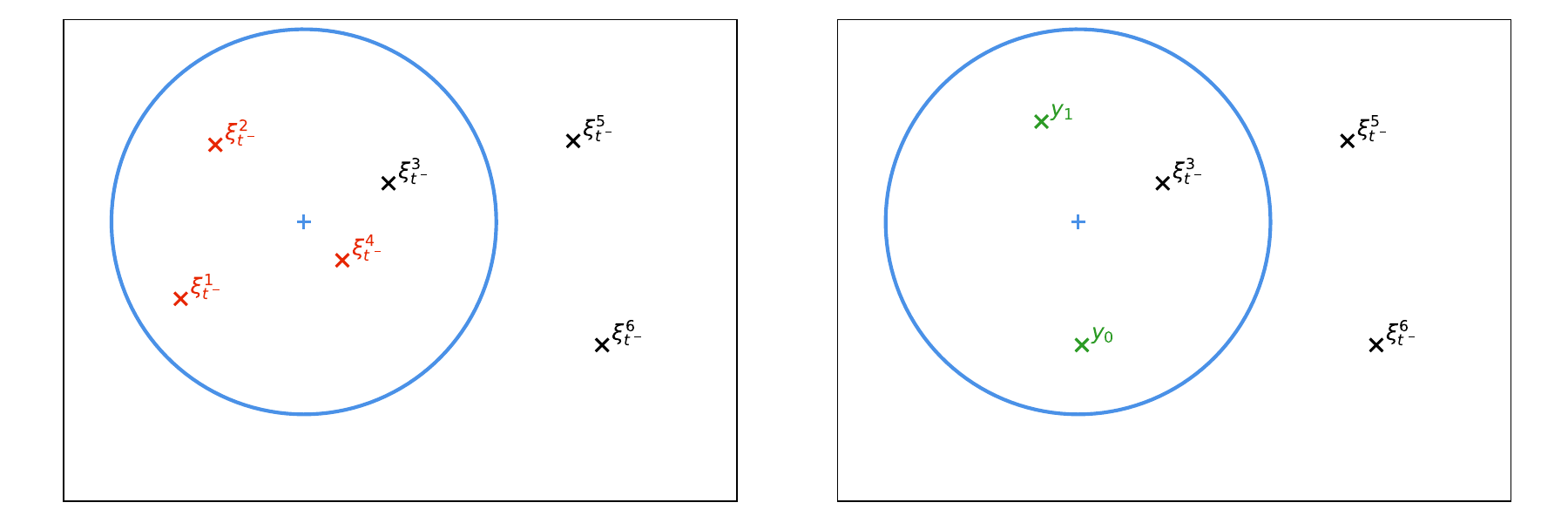}
    \includegraphics[width=0.95\textwidth]{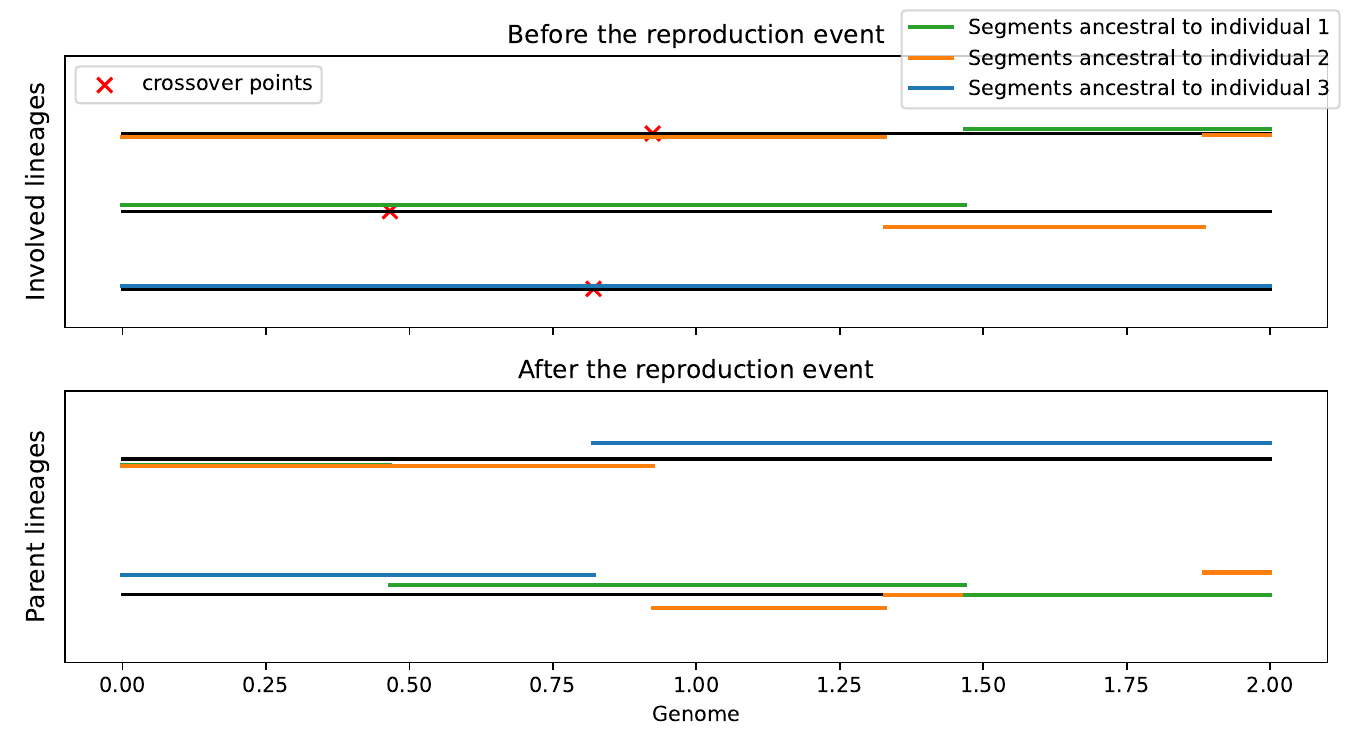}
    \caption{Illustration of a reproduction event in the genealogical process. a) Positions of the ancestors of the sampled individuals before (left) and after (right) a reproduction event. The lineages involved in the reproduction event are marked in red on the left, and the parent positions are marked in green on the right. b) For each involved lineage, we show the portions of genome that are ancestral to each sampled individuals carried by these lineages just before the reproduction event, with the crossover points marked by red crosses. c) Portions of genome of the parent lineages resulting from the reproduction event that are ancestral to the sampled individuals.}
    \label{fig:reproduction_event}
\end{figure}

An illustration of the state of $\mathcal{A}_t$ just before and just after a reproduction event is provided in Figure~\ref{fig:reproduction_event}.
Note that, with this definition, the rate at which a given lineage located at $\xi$ is involved in a new reproduction event is given by
\begin{equation*}
    \lambda(\xi) := \sum_{\alpha \in \PM} \frac{|\mathbb{H}_\alpha \cap B(\xi, r_\alpha)|}{V_{r_\alpha}},
\end{equation*}
where $|\cdot|$ denotes Lebesgue measure.
We then observe that, as soon as $|\xi^{(1)}| > R := \max(r_+,r_-)$, $\lambda(\xi) = 1$.
This justifies the scaling of the intensities of the Poisson random measures in Definition~\ref{def:ARG}, so that one time unit in this model corresponds to one generation.

\subsection{Segments of identity by descent}
We now formally define IBD segments in this model. Recall that an IBD segment is defined by comparing the genomes of a pair of sampled individuals. Without loss of generality, we can give the definition for the pair $(1,2)$, the definition is analogous for any pair. In words, an IBD segment is defined as a continuous portion of genome on which the two individuals share their most recent common ancestor unbroken by recombination events up to the time of this ancestor.

\begin{definition} \label{def:IBD_segment}
A genomic segment $[a,b[\, \subset [0,G]$ belongs to an IBD segment for the pair of samples $(1,2)$ if there exists $T_c > 0$ such that
\begin{enumerate}[i)]
    \item for all $0 \leq t < T_c$, there exists $1 \leq i \neq j \leq N_t$ such that
    \begin{equation*}
        1 \in \Block^i_t(\ell), \quad \text{ and } \quad 2 \in \Block^j_t(\ell), \quad \forall\, \ell \in [a,b[,
    \end{equation*}
    and
    \item there exists $1 \leq i \leq N_{T_c}$ such that, for all $\ell \in [a, b[$, $\lbrace 1, 2 \rbrace \in B^i_{T_c}(\ell)$.
\end{enumerate}
We then say that $[a,b[$ is an IBD segment if it belongs to an IBD segment and if it is maximal, i.e. any segment strictly containing $[a,b[$ does not belong to an IBD segment.
\end{definition}

Note that what is crucial in the above definition is that the time $T_c$ and the indices $i$ and $j$ in \textit{i)} are the same for all $\ell \in [a,b[$.
The time $T_c$ is then the coalescence time (or the time to the most recent common ancestor, TMRCA) for the IBD segment $[a,b[$. 
Note that every singleton $\lbrace \ell \rbrace$ satisfies the above definition (since the TMRCA is finite at every locus almost surely), so that IBD segments form a partition of the genome.
Adjoining IBD segments are thus delimited by recombination points that took place in one of the two lineages before the coalescence time.

According to a long-standing theoretical framework \cite{palamara_length_2012,ralph_geography_2013,palamara_inference_2013,carmi_renewal_2014,browning_accurate_2015,ni_probabilistic_2016,ringbauer_inferring_2017}, the expected number of IBD segments of a certain length can be derived as a function of the distribution of the coalescence time at a single locus. In particular, for $L \in (0, G)$, let $N_{\geq L}$ be the number of IBD segments of length at least $L$ found among an arbitrary pair of samples, and for an arbitrary locus $\ell \in [0, G - L]$ set
\begin{equation*}
    p_{L}(\ell) = \P{ [\ell, \ell +L[ \text{ belongs to an IBD segment}}.
\end{equation*}
The following then holds.
\begin{equation} \label{eq:expected_nb_segments}
    \E{N_{\geq L}} = - \deriv{}{L} \int_0^{G-L} p_L(\ell) d\ell.
\end{equation}

\begin{proof}
Let $(L_i, i \geq 1)$ denote the lengths of all the IBD segments for the pair of individuals.
Then we note that
\begin{equation*}
    \int_0^{G-L} \1{[\ell, \ell + L[ \text{ belongs to an IBD segment}} d\ell = \sum_{i \geq 1} \1{L_i \geq L} (L_i - L).
\end{equation*}
Taking expectations on both sides and differentiating with respect to $L$ yields \eqref{eq:expected_nb_segments}.
\end{proof}

This result is very useful because $p_L(\ell)$ can be computed using only the distribution of $(\mathcal{A}^\ell_t, t \geq 0)$, the genealogy of the sample at a single locus $\ell$.
Indeed,
\begin{equation} \label{eq:link_with_Tc}
    p_L(\ell) \approx \E{e^{-2L T_\ell}},
\end{equation}
where $T_\ell$ is the TMRCA at locus $\ell$ (see \cite{palamara_length_2012,carmi_renewal_2014,browning_accurate_2015,ringbauer_inferring_2017}).
In the model defined above, this formula is not exact, since recombination only takes place during reproduction events, hence $2 T_\ell$ should be replaced by the number of reproduction events involving one of the two lineages up to $T_\ell$ (since at each reproduction event the probability that no crossover point falls in $[\ell, \ell + L[$ is $e^{-L}$ for each lineage).
However, when $T_\ell$ is sufficiently large, the number of reproduction events involving each lineage will be close to $T_\ell$ (since each lineage is involved in a new reproduction event at rate 1), and the two formulas agree.

Note also that the distribution of $T_\ell$ does not depend on $\ell$, hence $p_L$ itself does not depend on $\ell$.
Plugging this in \eqref{eq:expected_nb_segments} yields
\begin{equation*} \label{eq:E_N_geq_L}
    \E{N_{\geq L}} = p_L - (G-L) \deriv{p_L}{L}.
\end{equation*}
In the following, we give a way to approximate $p_L$ under some conditions, but note already that, if we plug \eqref{eq:link_with_Tc} in the above, we obtain
\begin{equation*}
    \E{N_{\geq L}} \approx \E{e^{-2L T_0}} + 2(G-L) \E{T_0 e^{-2L T_0}},
\end{equation*}
as in \citep[equation (4)]{ringbauer_inferring_2017}.

\subsection{Wright-Malécot approximation with heterogeneous dispersal}

The main result of this section is an adaptation of the Wright-Malécot formula (see \cite{barton_inference_2013}) to our model, yielding an approximation of $p_L(\ell)$ under some assumptions.
To prove this formula we show that the rescaled trajectories of the lineages converge in distribution to continuous diffusions in $\R^2$. Here, due to the spatial heterogeneity in the radius of reproduction events, the resulting diffusion presents a singular drift at the interface $\lbrace x^{(1)} = 0 \rbrace$.

Let $(\xi_t, t \geq 0)$ denote the trajectory of the position of the genetic ancestor of an individual sampled at $\xi_0$ at one locus $\ell \in [0,G]$.
From Definition~\ref{def:ARG}, it can be seen that $(\xi_t, t \geq 0)$ is a jump Markov process with generator
\begin{equation} \label{def:generator_single_lineage}
    \mathcal{L}f(z) := \int_{\R^2} \Phi(z,y) \left( f(y) - f(z) \right) dy,
\end{equation}
where
\begin{equation} \label{def:Phi}
    \Phi(z,y) :=  \sum_{\alpha\in \PM} \frac{\abs{\mathbb{H}_\alpha \cap B(z,r_\alpha) \cap B(y,r_\alpha)}}{(V_{r_\alpha})^2}.
\end{equation}
Let us set $R := \max(r_-, r_+)$.
We then note that, when $\alpha z^{(1)} > R$, $B(z,r_{-\alpha}) \cap \mathbb{H}_{-\alpha} = \emptyset$, and $\mathcal{L}f(z)$ coincides with the generator of a symmetric random walk with variance
\begin{equation} \label{eq:variance}
    \sigma_\alpha^2 = \int_{\R^2} \frac{\abs{B(z,r_\alpha) \cap B(y,r_\alpha)}}{(V_{r_\alpha})^2} (y-z)^2 dy = \frac{r_\alpha^2}{2},
\end{equation}
(see \citep[Section~6.1]{barton_new_2010}).
In Appendix~\ref{app:cvg_to_skew_Bm}, we prove the following.

\begin{lemma} \label{lemma:cvg_single_lineage}
For $n \geq 1$, set
\begin{equation*}
    X^n_t := \frac{1}{\sqrt{n}} \xi_{nt}.
\end{equation*}
Assume that $X^n_0$ converges in distribution to $X_0$ taking values in $\R^2$.
Then $(X^n_t, t \geq 0)$ converges in distribution to $(X_t, t \geq 0)$ which solves
\begin{equation} \label{eds:skew_Bm}
    \left\lbrace
    \begin{aligned}
    &X^{(1)}_t = X^{(1)}_0 + \int_0^t \sigma(X_s) dB^{(1)}_s + \beta L^0_t(X^{(1)}), \\
    &X^{(2)}_t = X^{(2)}_0 + \int_0^t \sigma(X_s) dB^{(2)}_s,
    \end{aligned}
    \right.
\end{equation}
where $(B_t, t \geq 0)$ is standard Brownian motion taking values in $\R^2$,
\begin{equation*}
    \sigma(x) =  \sigma_\alpha \text{ for } x \in \mathbb{H}_\alpha, \qquad \beta = \frac{r_+^2 - r_-^2}{r_+^2 + r_-^2},
\end{equation*}
and $(L^0_t(X^{(1)}), t \geq 0)$ denotes the local time at 0 of $(X^{(1)}_t, t \geq 0)$ (see \citep[Chapter~VI]{revuz_continuous_2013}).
\end{lemma}

Observe that, in \eqref{eds:skew_Bm}, $t \mapsto L^0_t(X^{(1)}_t)$ increases only at the times for which $X^{(1)}_t = 0$, so that, in between those times, $(X_t, t \geq 0)$ behaves like Brownian motion with variance $\sigma_+^2$ or $\sigma_-^2$ depending on whether it is in $\mathbb{H}_+$ or $\mathbb{H}_-$.
The fact that \eqref{eds:skew_Bm} admits a unique solution for $\beta \in [-1, 1]$ is proved in \cite{lejay_constructions_2006} (Proposition~10) in the case $d=1$ and $\sigma_+ = \sigma_-$.
The extension to higher dimensions is straightforward and the case $\sigma_+ \neq \sigma_-$ can be treated via a time change depending on the first coordinate.
In \cite{portenko_diffusion_1979}, \cite{portenko_stochastic_1979}, it is proved that this process generates a Feller semigroup, which we denote by $(P_t, t \geq 0)$, i.e. for any continuous and bounded function $\phi : \R^2 \to \R$,
\begin{equation*}
    P_t\phi(x) = \E[x]{\phi(X_t)}.
\end{equation*}
It is shown in \cite{lejay_constructions_2006} that $(P_t \phi, t \geq 0)$ is also the solution to the initial value problem
\begin{equation} \label{pde_skew_diffusion}
    \left\lbrace
    \begin{aligned}
    & \partial_t \phi_t = \frac{\sigma_\alpha^2}{2} \Delta \phi_t,  \text{ in $\mathbb{H}_\alpha$}, \; \phi_t \in \mathcal{D}_\beta, \: t > 0, \\
    & \phi_0 = \phi,
    \end{aligned}
    \right.
\end{equation}
where $\mathcal{D}^\beta$ is the set of continuous functions $\phi : \R^2 \to \R$, twice continuously differentiable on each half-space, such that
\begin{equation*}
    (1+\beta) \deriv*{\phi}{x_1}{0^+} = (1-\beta) \deriv*{\phi}{x_1}{0^-},
\end{equation*}
with all the other derivatives continuous on $\R^2$.

We are now ready to state our main result.
Let us now assume that, for $\alpha \in \PM$,
\begin{equation*}
    u_\alpha = \frac{\upsilon_\alpha}{N},
\end{equation*}
for some $\upsilon_\alpha \in (0, 1]$ and $N \geq 1$.
Take $\phi, \psi : \R^2 \to \R$ two probability density functions on $\R^2$, and set, for $d > 0$,
\begin{equation} \label{def:phi_d}
    \phi_d(x) := \frac{1}{d^2} \phi\left( \frac{x}{d} \right), \quad \psi_d(x) := \frac{1}{d^2} \psi\left( \frac{x}{d} \right).
\end{equation}
Then, for $d > 0$ and $L > 0$, let $p_{L, d}(\phi, \psi)$ denote the probability that $[\ell, \ell + L[$ belongs to an IBD segment for a pair of individuals sampled from the probability density functions $\phi_d$ and $\psi_d$.

\begin{theorem} \label{thm:WM_formula}
If $d \gg 1$, $N \gg \log(d)$ and $d^2 L \asymp 1$, then
\begin{equation} \label{eq:WM_formula}
    N p_{L, d}(\phi, \psi) \approx \int_0^\infty \int_{\R^2} \frac{1}{\mathcal{N}(x)} e^{-2 d^2 L t} P_t\phi(x) P_t\psi(x) dx dt,
\end{equation}
where $\mathcal{N} : \R^2 \to \R_+$ is defined by
\begin{equation} \label{def:N}
    \mathcal{N}(x) := \frac{2}{\upsilon_\alpha V_{r_\alpha}} \text{ for } x \in \mathbb{H}_\alpha,
\end{equation}
and $(P_t, t \geq 0)$ is defined as above.
\end{theorem}

Here, $d^2 L \asymp 1$ means that $d^2 L$ should be neither too large nor too small, and $\approx$ actually means that the left-hand side converges along sequences $N_n$, $L_n$, $d_n$ such that $N_n, d_n \to \infty$, $\frac{\log(d_n)}{N_n} \to 0$ and $d_n^2 L_n \to \gamma$ for some $\gamma \in (0, \infty)$, and that the limit is given by the obvious limit of the right-hand side (see Appendix~\ref{app:Wright_Malecot}).

Theorem~\ref{thm:WM_formula} tells us that the expected number of IBD segments longer than $L$ shared by a pair of individuals depends on a small set of demographic parameters, namely the dispersal variance $\sigma_\alpha^2$ and the so-called neighbourhood size $\mathcal{N}_\alpha := 2N / (\upsilon_\alpha V_{r_\alpha}) = 2 / (u_\alpha V_{r_\alpha})$ on each side of the interface.
It can be seen from the proof of this result that \eqref{eq:WM_formula} holds under more general models, provided one replaces the values of these parameters by their equivalent in the new model (e.g. in a stepping stone model as in \cite{nagylaki_clines_1976}).

\begin{figure}
    \centering
    \makebox[\textwidth][c]{\includegraphics[height=0.4\textheight]{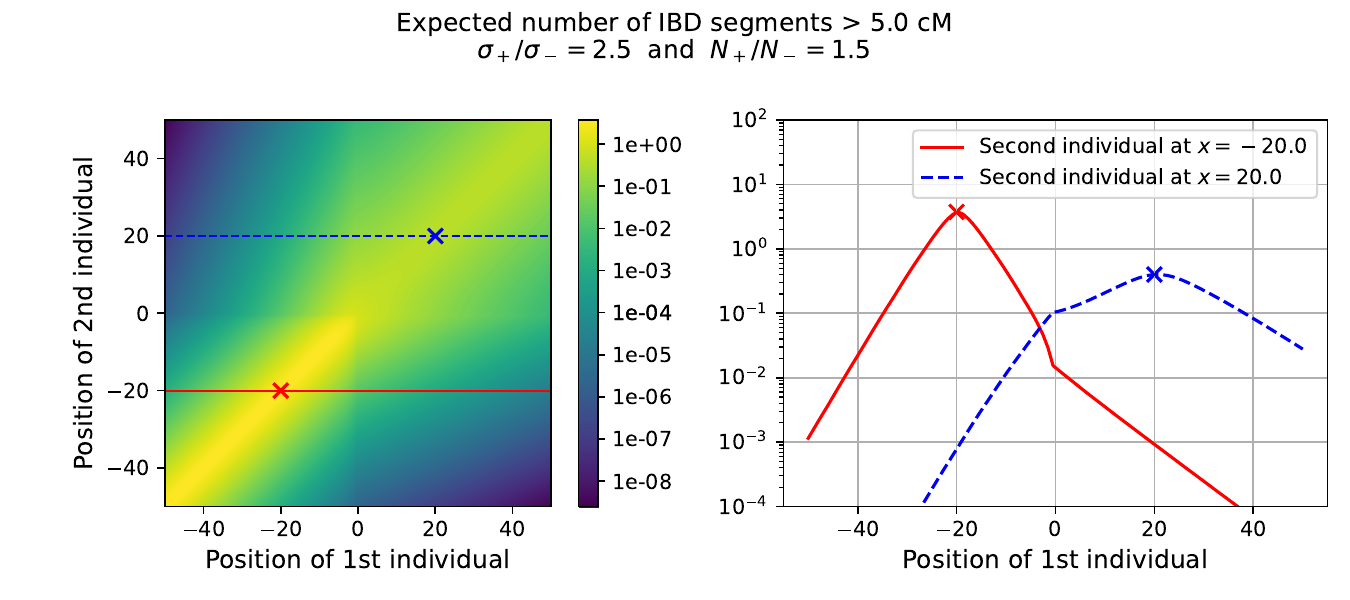}}
    \caption{Expected number of IBD segments longer than 5cM shared by a pair of individuals sampled from a linear transect perpendicular to the interface between $\mathbb{H}_+$ and $\mathbb{H}_-$, according to \eqref{eq:WM_E_N_geq_L}. The left graphic shows a heatmap of the expected number of segments for all pairs of individuals along this transept, with the position of the first individual along the transept given in the $x$ axis, and the position of the second individual in the $y$-axis.
    The right graph shows the expected number of IBD segments ($y$-axis) for two fixed positions of the second individual (represented by the dashed lines on the left graphic), with the position of the first individual along the transept in the $x$-axis.
    Parameter values for this plot are $ \sigma_- = 1$, $ \sigma_+ = 2.5$, $\mathcal{N}_- = 10$, $ \mathcal{N}_+ = 15$, and a genome length of 50M. We can see a sharper decrease of the expected number of IBD segments with the distance between the individuals in $\mathbb{H}_{-}$, where dispersal is weaker, but also a larger expected number of segments for nearby pairs, as the neighbourhood size is also smaller. The right graphic also shows the difference in the exponential decay of expected IBD sharing between the two half-spaces, with a non-trivial transition at the interface. In particular, we observe a break in the slope of the two curves at the interface, as expected from \eqref{eq:WM_E_N_geq_L} and \eqref{pde_skew_diffusion}.}
    \label{fig:expected_ibd}
\end{figure}

Plugging \eqref{eq:WM_formula} in \eqref{eq:E_N_geq_L}, we obtain
\begin{equation} \label{eq:WM_E_N_geq_L}
    \E{N_{\geq L}} \approx \int_0^\infty \int_{\R^2} \frac{1}{N \mathcal{N}(x)} \left[ 1 + 2 (G-L) d^2 t \right] e^{-2d^2 L t} P_t\phi(x) P_t\psi(x) dx dt.
\end{equation}
We use this formula to estimate the expected number of IBD segments whose length falls in a certain interval in our inference method (see below).
The values taken by this function are plotted on Figure~\ref{fig:expected_ibd} for a particular set of pairs of individuals (sampled along a linear transept crossing the interface between the two half-spaces).

Recall that, in the above formula, $\phi$ and $\psi$ are two probability density functions, which encode the distribution of the sampling location of the two individuals.
They thus account for the uncertainty in the actual location of the samples.
If the sampling location is known exactly (as is the case in our simulations below), we can instead replace $P_t\phi(x)$ and $P_t\psi(x)$ in \eqref{eq:WM_E_N_geq_L} by $p_t(x, x_1)$ and $p_t(x, x_2)$, where $x_1$ and $x_2$ are the two sampling locations and $p_t(x,y)$ is the transition density of the diffusion defined in \eqref{eds:skew_Bm}.

\section{Inference method} \label{sec:method}

Our inference scheme aims to estimate the demographic parameters of the model, $\theta = \lbrace (\sigma_\alpha^2, \mathcal{N}_\alpha), \alpha \in \PM \rbrace$, using the observed IBD segments between pairs of sampled individuals.
As noted in \cite{ringbauer_inferring_2017}, the joint distribution of the IBD segment lengths in the sampled individuals is generally intractable since it involves the law of the genealogy of the samples along the whole genome (in other words, the law of $(\mathcal{A}^\ell, \ell \in [0,G])$). 
As a result, we cannot hope to compute the exact likelihood of a set of parameters $\theta$. 
We thus adapt to our setting the composite likelihood approach proposed in \cite{ringbauer_inferring_2017}.

Let us fix a disjoint family of segment length intervals $\lbrace [L_k, L_k + \Delta L_k[, k \in \mathcal{K} \rbrace$, and let $N_{i,j,k}$ denote the number of IBD segments with length in the interval $[L_k, L_k + \Delta L_k[$ found in the pair of individuals $(i,j)$, for $i \neq j$ (in practice, we shall take a fixed value for $\Delta L_k$ and $L_{k+1} = L_{k} + \Delta L$).
Also let $\lbrace x_1, \ldots, x_n \rbrace$ be the locations of the $n$ sampled individuals.

To approximate the likelihood of the parameters, we assume that the numbers of IBD segments in each length interval found in any given pair of samples are independent and follow Poisson distributions with parameters
\begin{equation*}
    \lambda_{i,j,k}(\theta) := \Etheta[(x_i, x_j)]{N_{\geq L_k}} - \Etheta[(x_i,x_j)]{N_{\geq L_k + \Delta L_k}}, \quad 1 \leq i \neq j \leq n, \quad k \in \mathcal{K},
\end{equation*}
where $\Etheta[(x,y)]{\cdot}$ denotes the expectation under the model with parameters given by $\theta$, with two lineages started at locations $x$ and $y$.
The composite likelihood of a set of parameters $\theta$ given the observations
\begin{equation*}
    X := \lbrace N_{i,j, k}, k \in \mathcal{K}, 1 \leq i \neq j \leq n \rbrace
\end{equation*}
is then given by
\begin{equation*}
    cL(\theta, X) := C \prod_{k \in \mathcal{K}} \prod_{1 \leq i < j \leq n} \lambda_{i,j,k}(\theta)^{N_{i,j, k}} e^{-\lambda_{i,j,k}(\theta)},
\end{equation*}
where $C > 0$ is a normalising constant which does not depend on $\theta$.
The maximum likelihood estimator of the parameters, denoted by $\hat{\theta}$, is then obtained as
\begin{equation*}
    \hat{\theta} = \argmax_{\theta} cL(\theta, X).
\end{equation*}
We estimate the variance-covariance matrix of $\hat{\theta}$ using the inverse of the observed Fisher information matrix, and derive 95\% confidence intervals accordingly, even though these will not be true 95\% confidence intervals since the Fisher information matrix derived from a composite likelihood underestimates the variance of the estimators \cite{lindsay_composite_1988}. 

It is interesting to gather the observed IBD segments in discrete length intervals as it fixes the values of $L$ at which the function \eqref{eq:WM_E_N_geq_L} has to be evaluated, irrespective of the sample size. If several individuals have been sampled at approximately the same location, it is possible to reduce the time needed to evaluate the composite likelihood by summing all the IBD segments of each length interval involving lineages originating from the same pairs of locations, and by weighting the quantities $\lambda_{i,j,k}$ by the number of pairs of individuals between the different sampling locations.

The assumption that the numbers of observed segments in each length interval are independent is equivalent to assuming that the total number of segments in each pair of individuals is a Poisson random variable, and that the lengths of the different segments are i.i.d.
In this approximation, we lose all the correlations in the genealogies along the genome of the individuals, which are induced by the fact that different loci share a common history in the recent past.
Nonetheless, if successive IBD segments are sufficiently far apart along the genome, then it only takes a few generations for recombination to take place between them in both lineages and so they likely coalesce through different paths, therefore correlations between their lengths should be small.

The composite likelihood also does not keep track of correlations between different pairs of individuals, for example, two individuals who already share IBD segments with a third individual are often more likely to share IBD segments between them than an average pair of individuals.
Again, if IBD segments are sufficiently rare, the probability that several segments come from the same ancestry should remain small.

In any case, as already noted in \cite{ringbauer_inferring_2017}, provided the approximation of the marginals by Poisson random variables is accurate, maximizing the composite likelihood  provides consistent and asymptotically normal estimates of the model parameters \cite{lindsay_composite_1988}.
One issue is that the confidence intervals obtained from the curvature of the composite log-likelihood surface at its maximum will be too tight due to the assumption of independence of IBD blocks \cite{coffman_computationally_2016}, see Section~\ref{sec:results}.

In order to compute $\lambda_{i,j,k}(\theta)$, we need to evaluate \eqref{eq:WM_E_N_geq_L}, and hence to solve \eqref{pde_skew_diffusion}. Contrary to the usual diffusion equation, this equation does not admit simple analytical solutions. We thus use a numerical approximation of $P_t\phi$.
To obtain this approximation, we compute the distribution of a random walk on a fine discrete grid $(Y_t, t \geq 0)$, chosen so that this random walk admits the skew diffusion \eqref{eds:skew_Bm} as a scaling limit, and we approximate
\begin{equation*}
    P_t \phi(x) = \E[x]{\phi(X_t)} \approx \E[\mu_x]{\phi(Y_t)},
\end{equation*}
where $\mu_x$ is a probability distribution supported on the four grid points closest to $x$. Details of this approximation are given in Section~\ref{sec:random_walk} of the supplementary material.

\section{Results on simulated IBD segment sharing data} \label{sec:results}
We tested our inference method on simulated IBD segments using two simulation models: the spatial $\Lambda$-Fleming-Viot process (see Definition~\ref{def:ARG}) and a stepping stone model on a grid with heterogeneous deme size and migration.

The stepping stone model consists of a two-dimensional finite grid of demes. 
Each deme contains a fixed number of chromosomes. At each generation, each chromosome descends from a random pair of chromosomes in the previous generation, either from the same deme or a neighbouring deme, chosen with probabilities according to a fixed migration kernel. Recombination between the two parental chromosomes occurs according to the model described in Section~\ref{sec:model}. We chose a migration kernel similar to the one studied in \cite{nagylaki_clines_1976} (but for a two-dimensional grid, see Figure~\ref{fig:kernel}).

In each simulation run, 25 sampling locations were chosen, forming a regular grid around the origin, and 10 lineages were started at each of these locations, each carrying a genome of 25M (for the simulations using the spatial $\Lambda$-Fleming-Viot process) or 50M (for the simulations on the discrete grid).
We then simulated the ancestral recombination graph according to either model for a fixed number of generations (chosen large enough that long IBD segments are extremely unlikely to coalesce beyond the simulation time range). At each coalescence event, we recorded all resulting IBD segments longer than 4cM.

\begin{table}[ht]
    \centering
    \begin{tabular}{l|rrrr}
        Scenario & 1 & 2 & 3 & 4 \\
        \hline
        $\sigma_-$ & 0.5 & 0.5 & 0.4 & 0.4 \\
        $\sigma_+$ & 0.5 & 0.5 & 0.8 & 0.8 \\
        $\mathcal{N}_-$ & 80 & 40 & 80 & 40 \\
        $\mathcal{N}_+$ & 80 & 80 & 80 & 80
    \end{tabular}
    \caption{Parameters used for the simulations. We note that scenario~1 corresponds to the homogeneous model previously studied in \cite{ringbauer_inferring_2017}. Scenario~2 presents homogeneous dispersal but heterogeneous effective neighbourhood size, while Scenario~3 is the reverse. Finally, Scenario~4 presents both heterogeneous dispersal and effective neighbourhood size.}
    \label{table:simu_params}
\end{table}

\begin{figure}[ht]
    \centering
    \includegraphics[width=\textwidth]{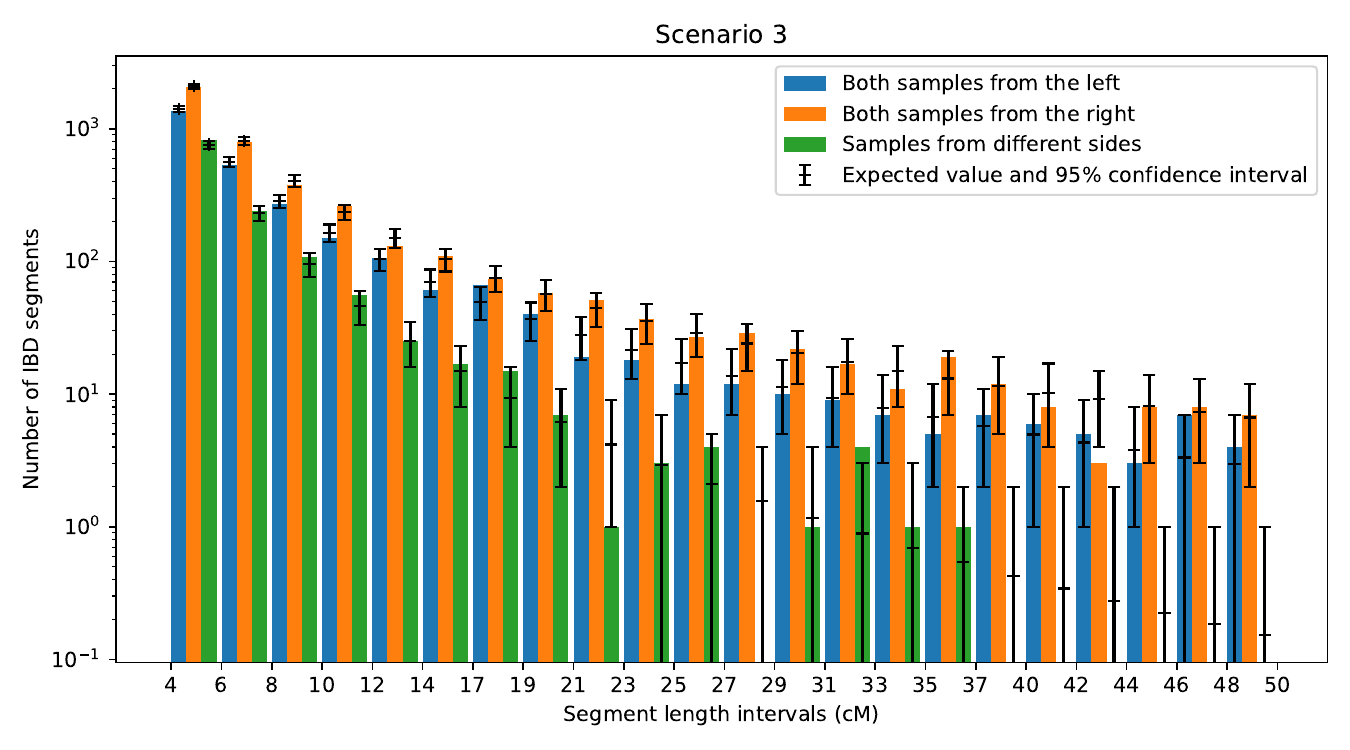}
    \caption{Number of IBD segments of different lengths observed in one simulation of the SLFV for scenario~3 (see Table~\ref{table:simu_params}).
    The bar plot shows the sum of all observed IBD segments of each length interval for all pairs of individuals respectively both on the left side of the interface, both on the right side, and on different sides. The expected value computed with \eqref{eq:WM_E_N_geq_L} with the simulation parameters and the associated 95\% confidence interval of the Poisson distribution are shown in black for each length interval and each sampling configuration. The same graphs for scenarios 1, 2 and 4 are provided in Figure~\ref{supp:fig:observed_ibd} of the supplementary material.}
    \label{fig:observed_ibd}
\end{figure}

We chose four main sets of parameters consisting of dispersal variance and effective neighbourhood size on each side of the interface, given in Table~\ref{table:simu_params} (additional scenarios were run and are presented in the supplemental material, see Figure~\ref{fig:simu_results_varying_pop_size}).
For each set of parameters (or each scenario, as we call them), we performed 20 independent simulations.
Figure~\ref{fig:observed_ibd} displays the numbers of IBD segments in one simulation, along with the expected values computed with \eqref{eq:WM_E_N_geq_L}.

\begin{figure}[tp]
    \centering
    \includegraphics[width=0.95\textwidth]{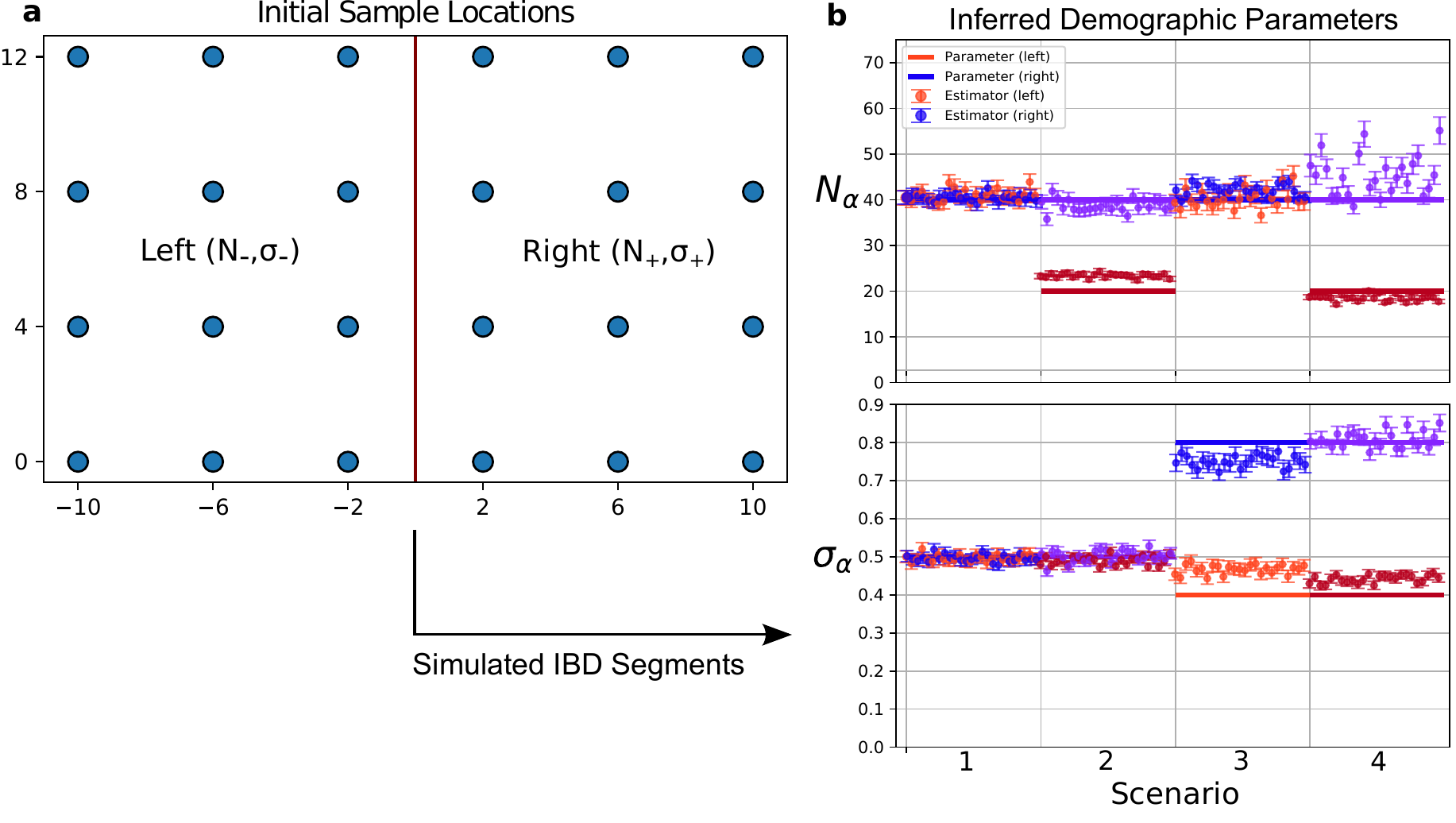}
    \caption{\textbf{a}: Sampling locations on a 2D spatial grid. The vertical red line depicts the interface between $\mathbb{H}_-$ and $\mathbb{H}_+$. For each replicate run of a scenario, ten genomes are sampled at each initial sampling location, and IBD segments between those samples are simulated using the stepping stone grid simulation scheme as described in the main text.
  \textbf{b}: We depict estimated (dots) vs. simulated (horizontal lines) parameters in 25 replicates of four different scenarios for the stepping stone model.
    Confidence intervals computed using the Fisher information matrix are shown by the error bars around each estimate. The top panel depicts the effective neighbourhood size and the bottom panel the dispersal parameters. The parameters on the left ($\mathcal{N}_-$ and $\sigma_-$) of the interface are depicted in red and those on the right ($\mathcal{N}_+$ and $\sigma_+$) in blue. In the above graph, $N_\alpha$ denotes the diploid effective population size, which is half the haploid effective population size, hence the factor 2 between the values above and those given in Table~\ref{table:simu_params}.}
    \label{fig:simu_results_stepping_stone_fixed_beta}
\end{figure}

\begin{figure}[tp]
    \centering
    \includegraphics[width=0.7\textwidth]{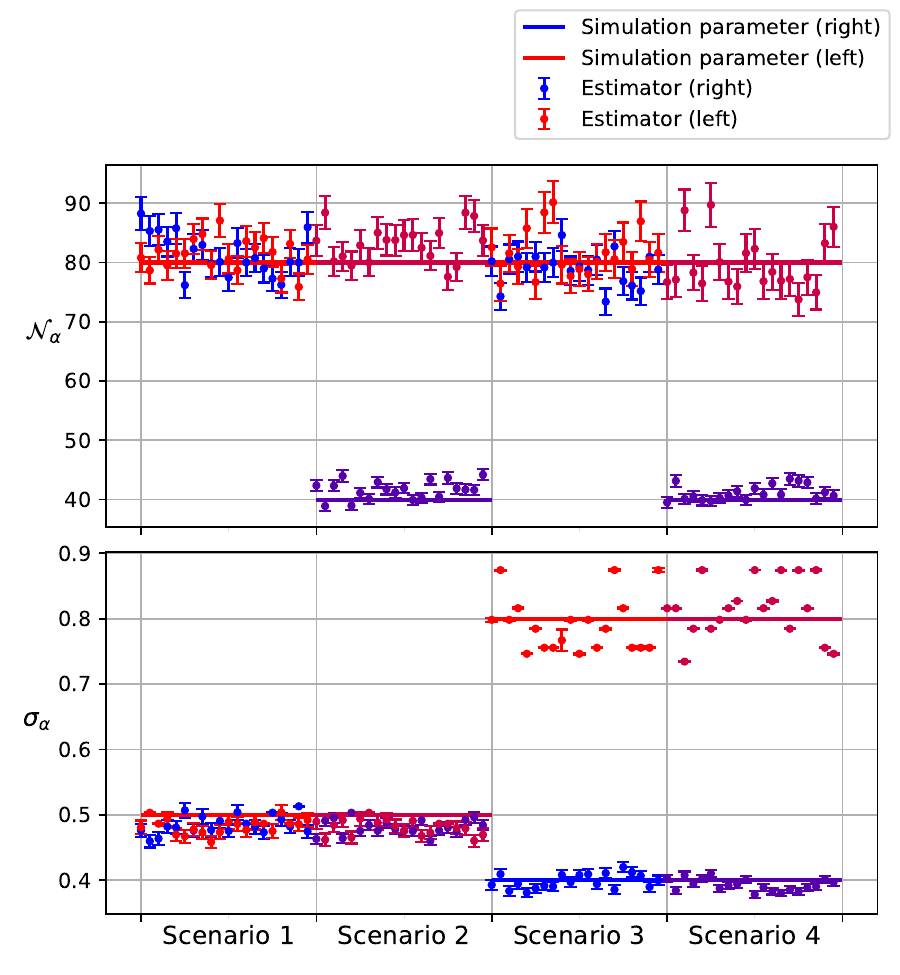}
    \caption{Simulation results using the spatial $\Lambda$-Fleming-Viot process of Definition~\ref{def:ARG}, using the simulation parameters  of Table~\ref{table:simu_params}. These simulations were run with $G$ = 25M (instead of 50M in the grid model) to limit the simulation time. Each scenario was independently simulated 20 times, and the same inference procedure was applied to the resulting IBD segments. Simulated and estimated parameters are shown as in Figure~\ref{fig:simu_results_stepping_stone_fixed_beta}.}
    \label{fig:results_slfv}
\end{figure}

We then applied the inference method to each replicate and compared the values of the parameters used in each scenario to the estimators and confidence intervals obtained with our method (see Figure~\ref{fig:simu_results_stepping_stone_fixed_beta} and~\ref{fig:results_slfv}). We found that our method performs well in a variety of settings. In scenario~1, where the parameters are the same in both regions, the estimates for both sites largely overlap and agree with the parameters used in the simulations. 
The estimator is also able to robustly distinguish between heterogeneous neighbourhood size (by a factor of 2, scenario 2) and heterogeneous dispersal coefficients (by a factor of 2, scenario 3), and is also able to detect cases where both types of parameters are heterogeneous (scenario 4). These observations show that our inference method is able to detect sharp variations in either or both kinds of parameters (here by a factor of 2) even without a very dense sampling scheme.

The estimates appear to be slightly biased in some cases, tending to underestimate the difference between the parameters on each side of the interface (see for instance the results of scenario 3 for the stepping stone model).
Recall, on the other hand, that since we use a composite likelihood, the confidence intervals estimated from the Fisher information matrix are tighter than they should be, so it is not surprising that the values used in the simulations sometimes lie outside this interval.

Assuming constant population density, even over the span of a few dozen generations, can be problematic, especially when analysing data from human samples.
However, the inference method developed in \cite{ringbauer_inferring_2017} can also accommodate some predefined forms of temporal variations of the effective neighbourhood size.
Indeed, the only modification needed in the above-described scheme is to replace $\mathcal{N}(x)$ by some function of time and space $\mathcal{N}(t,x)$ in \eqref{eq:WM_E_N_geq_L}, where $\mathcal{N}(t,x)$ takes a simple form with a small number of parameters, and to infer these parameters in the inference step. 

Note that, in \eqref{eq:WM_E_N_geq_L}, time runs backwards and has been scaled by $d^2$ (see Appendix~\ref{app:Wright_Malecot}), so that $\mathcal{N}(t,x)$ corresponds to the effective neighbourhood size $t / d^2$ generations before present.
As an illustration, we assumed that $\mathcal{N}(t,x)$ was of the form $\mathcal{N}_\alpha t^{-\gamma}$ in each half-space $\mathbb{H}_\alpha$ (with the same $\gamma$ in both half-spaces, and we ran simulations of the discrete grid model with this varying population size.
The results of the inference procedure on these simulated data is presented in Figure~\ref{fig:simu_results_varying_pop_size} of the supplementary material, for nine different sets of parameters.
In every tested scenario, the addition of the growth parameter did not greatly affect the inference results, and all five parameters were successfully inferred. 

\section{Discussion} \label{sec:discussion}
When analyzing genetic data issued from heterogeneous populations, it can be tempting to divide the samples into smaller subsamples, corresponding to spatial regions one expects to be more homogeneous, and to analyse each smaller data set independently of one another, ignoring the interactions between the subsamples and the potential border effects near the interface between the different regions. Doing so, however, can result in a substantial loss of information and of power, and might introduce biases.

Another promising approach is to consider that the demographic parameters vary continuously in the considered spatial domain and to estimate these parameters at each point of some dense spatial grid.
This approach was adopted in \cite{al-asadi_estimating_2019}, also using as signal segments of shared haplotypes, both on simulated data and on the POPRES data set (human SNP array data from across Europe \cite{nelson_population_2008}).
The resulting migration and population size surfaces indeed displayed significant spatial (and even temporal) heterogeneity.
Such methods can however be sensitive to uncertainties regarding the sampling locations, and might also require sufficiently regular sampling to obtain consistent estimates.
In addition, due to a large number of inferred parameters, \cite{al-asadi_estimating_2019} need to introduce a regularising penalisation in their estimation to avoid overfitting. 
Such a penalization is hard to balance correctly, as it can reduce the ability to detect sharp changes in the effective demographic parameters.

We thus propose a complementary method which accounts for spatial heterogeneity while exploiting nearly all the information from IBD blocks in the sampled genotypes. As for any other demographic inference method, the values of the estimated parameters should be seen as a concise summary statistic of a very complex history of migration and demography. Direct interpretations in terms of actual distances travelled by individuals require caution (see for example the caveats described in \cite{bradburd_spatial_2019}).
Our simulations show, however, that the proposed method can reliably detect discontinuities and qualitative variations in migration patterns and effective population density.
Moreover, as expression \eqref{eq:WM_E_N_geq_L} does not depend on the details of the stochastic model, the method is robust to some amount of model misspecification.
We were for example able to use it to infer the parameters of a stepping-stone model on a regular grid.
In addition, the computation time needed to infer one set of parameters is very short, of the order of only a few minutes for a sample of about 250 individuals.

When analysing real genomes, IBD segments cannot be directly observed, and one has to infer the lengths of IBD segments from dense single nucleotide or whole genome sequences \cite{browning_fast_2011,browning_probabilistic_2020,zhou_fast_2020}.
This estimation is prone to errors, including false positives (calling IBD segments in the absence of recent shared ancestry), false negatives (some segments might go undetected due to recent mutations) and errors regarding the exact segment length.
Ralph and Coop \cite{ralph_geography_2013} developed a model to estimate the rate of these detection errors, by running the IBD detection algorithm on modified data sets with ``spiked-in'' IBD segments.
Recently, \cite{browning_probabilistic_2020} introduced a Bayesian framework to detect long shared haplotype segments, whose output is the posterior distribution of the segment's endpoints given the observed sequences. These detection errors can be accounted for in the composite likelihood function in a straightforward way, by modifying the expectations $\lambda_{i,j,k}(\theta)$ to expectations that include the error model (see e.g. \cite{ralph_geography_2013,ringbauer_inferring_2017}).

In addition, recombination events that would end IBD segments might be impossible to detect if the two resulting lineages rapidly coalesce with each other just after recombination. Such recombination events are called ineffective in \cite{barton_inference_2013}, and increase the length of apparent IBD segments. Estimating the rate of such events is difficult, and might depend in a non-trivial way on the details of the model. Nonetheless, if the neighbourhood size is large enough, the probability that the two lineages coalesce with each other very quickly becomes very small, so the effects of ineffective recombination become negligible for all but extreme demographic scenarios.
Moreover, two or more short IBD segments that are close to each other can be inferred as a single longer IBD segment (because too few mutations take place between them), potentially biasing the IBD segment length distribution.
However \cite{chiang_conflation_2016} showed that this is effect becomes negligible for segments longer than 2cM.

In many empirical datasets such as the POPRES data set, the only geographical information on the location of the samples is the country of origin or other geographical unit, resulting in substantial uncertainty concerning the precise sampling locations. This uncertainty can be accounted for when computing the expected IBD sharing via the sampling distributions $\phi$ and $\psi$ in \eqref{eq:WM_E_N_geq_L}.
Ideally, one would even compute the composite likelihood for each set of possible sampling locations, and then sum over the product of the sampling distributions, to account for correlations in IBD sharing induced by the sampling locations. However, this would significantly increase the computation requirements of the inference method. Moreover, if the uncertainty on the sampling locations is of the same order as $\sigma_\alpha$, the variance of the position of lineages due to this uncertainty will be negligible after a few generations.

Building upon the approach introduced here, many extensions can be considered, such as non-straight interfaces or including more than two regions with different demographic parameters. For example, it would be possible to include the physical limits of the spatial domain in which the population evolves by replacing $P_t$ in \eqref{eq:WM_E_N_geq_L} by the semigroup of a similar diffusion reflected at the boundary of this domain.

Our results thus pave the way for a flexible inference framework for spatially heterogeneous populations, using the growing set of haplotype sharing data. This method and its proposed extensions could be used to test different demographic scenarios, to estimate the regions where demographic parameters change abruptly, etc. The combination of contemporary genetic data with ancient DNA is another exciting direction which offers many potential applications, including the investigation of changes in the patterns of dispersal or demographic heterogeneity.

\section*{Code availability}

The simulations and analysis described in this manuscript were implemented using the Python programming language. The code to simulate the discrete grid model and to perform the IBD-based inference is available at \url{https://github.com/hringbauer/IBD-Analysis}. The packages used to simulate the coalescent model in continuous space are available at \url{https://github.com/rforien/SLFV}.

\section*{Competing interests statement}
The authors declare that they have no competing interest relating to the content of this manuscript.

\section*{Acknowledgements}
The authors wish to thank Nick H. Barton, Alison Etheridge and Amandine V\'eber for their guidance during the early stages of this project. 
Both R.F. and H.R. are grateful to the Coop Lab for hosting them at different times and for many fruitful discussions throughout these periods.
R.F. is partially supported by the chaire Modélisation Mathématiques et Biodiversité VEOLIA-Ecole Polytechnique-MnHn-FX. HR is supported by funding from the Max Planck Society. Funding was also provided by the National Institutes of Health (NIH R35 GM136290 awarded to GC).

\section*{Author Contributions}
We annotate author contributions using the CRediT Taxonomy labels (\url{https://casrai.org/credit/}). 

\begin{itemize}
    \item Conceptualization (Design of study) -- RF
    \item Software -- RF, HR
    \item Formal Analysis -- RF
    \item Visualization -- RF, HR
    \item Writing (original draft preparation) -- RF
    \item Writing (review and editing) -- All authors
    \item Supervision -- GC 
\end{itemize}

\appendix
\renewcommand{\thesection}{\Alph{section}}
\titleformat{\section}{\bf}{Appendix \thesection.}{1em}{}

\section{Wright-Mal\'ecot formula} \label{app:Wright_Malecot}

Here we provide the proof for Theorem~\ref{thm:WM_formula}.
More precisely, let $(N_n, n \in \N)$, $(d_n, n \in \N)$ and $(L_n, n \in \N)$ be positive sequences such that, as $n \to \infty$,
\begin{equation} \label{parameter_scaling}
    d_n \to \infty, \quad \frac{N_n}{\log(d_n)} \to \infty, \quad d_n^2 L_n \to \gamma,
\end{equation}
for some $\gamma \in (0, \infty)$.
Also assume that, if $(\xi_t, t \geq 0)$ is a random walk with generator given by \eqref{def:generator_single_lineage}, then for any $\tilde{R} > 0$, there exists $ C > 0 $ such that
\begin{equation} \label{local_clt_bound}
    \sup_{x \in \R^2} \P[\xi_0]{\xi_t \in B(x, \tilde{R})} \leq \frac{C}{t}, \quad \forall\, t > 0, \xi_0 \in \R^2.
\end{equation}
(This bound holds for simple random walks by the local central limit theorem \citep[Chapter~2]{lawler_random_2010}.)
We then consider a sequences of processes $\lbrace (\mathcal{A}^n_t, t \geq 0), n \in \N \rbrace$ defined as in Definition~\ref{def:ARG}, with
\begin{equation*}
    u_\alpha^n = \frac{\upsilon_\alpha}{N_n},
\end{equation*}
and we let $p^n_{L_n, d_n}(\phi, \psi)$ denote the probability that $[\ell, \ell + L_n[$ belongs to an IBD segment, starting from two lineages distributed according to $\phi_{d_n}$, $\psi_{d_n}$ (defined in \eqref{def:phi_d}).
We then aim to prove that, as $n \to \infty$,
\begin{equation} \label{cvg_pn}
    N_n p^n_{L_n, d_n}(\phi, \psi) \to \int_0^\infty \int_{\R^2} \frac{1}{\mathcal{N}(x)} e^{-2 \gamma t} P_t\phi(x) P_t\psi(x) dx dt,
\end{equation}
where $\mathcal{N}$ and $P_t$ are defined in the statement of Theorem~\ref{thm:WM_formula}.
Note that, without loss of generality, we can assume that $d_n = \sqrt{n}$.
The bound \eqref{local_clt_bound} is used in the proof of Lemma~\ref{lemma:cvg_lineages} below, to control the expectation of the time that two lineages spend close to each other.

For $ i \in \lbrace 1, 2 \rbrace $, let $\xi^{n,i}_t$ denote the position at time $t$ of the genetic ancestor of the $i$-th individual at locus $\ell$.
Let $T_c^n$ be the time at which these two lineages coalesce, and let $T^n_{rec}$ be the first time at which a recombination event takes place in either lineage on the segment $]\ell, \ell + L_n[$.
Up to $T^n_{rec}$, the two lineages at $\xi^{n,1}_t$ and $\xi^{n,2}_t$ are also genetic ancestors of the sampled individuals for the whole segment $[\ell, \ell + L_n[$.
Thus, $p^n_{L_n, d_n}(\phi, \psi)$ is also the probability that $T^n_c < T^n_{rec}$.

Given that a lineage is involved in a reproduction event, the probability that the segment $[\ell, \ell + L_n[$ does not undergo recombination is $e^{-L_n}$.
Moreover, given that the two lineages are involved in the same reproduction event, and that neither lineage sees a recombination in $[\ell, \ell + L_n[$, they coalesce with probability $1/2$ (i.e. they inherit this segment from the same parent with probability $1/2$).
For $z_1 \in \R^2$, $z_2 \in \R^2$, let $q^n( z_1,z_2 )$ denote the probability that $T^n_c < T^n_{rec}$, starting from two lineages in $z_1$ and $z_2$.
We can then write, using \eqref{def:phi_d},
\begin{equation} \label{eq_pn_qn}
    p^n_{L_n, d_n}(\phi, \psi) = \int_{\R^2 \times \R^2} q^n(\sqrt{n} x_1, \sqrt{n} x_2) \phi(x_1) \psi(x_2) dx_1 dx_2.
\end{equation}
In addition, $q^n$ is symmetric and, by the Markov property, it satisfies the following equation
\begin{multline*}
    \sum_{1 \leq i \neq j \leq 2} \sum_{\alpha \in \PM} \frac{1}{u^n_\alpha V_{r_\alpha}} \int_{\mathbb{H}_\alpha \cap B(z_i, r_\alpha)} u^n_\alpha \left[ 1 - u^n_\alpha \1{B(z_j, r_\alpha)}(x) \right] \\ \times \left\lbrace e^{-L_n} \frac{1}{V_{r_\alpha}} \int_{B(x, r_\alpha)} \left( q^n(z_j, y) - q^n(z_1, z_2 ) \right) dy - (1-e^{-L_n}) q^n( z_1, z_2 ) \right\rbrace dx \\
    + \sum_{\alpha \in \PM} \frac{(u^n_\alpha)^2}{u^n_\alpha V_{r_\alpha}} \int_{\mathbb{H}_\alpha \cap B(z_1, r_\alpha) \cap B(z_1, r_\alpha)} \Bigg\lbrace \frac{1}{2} e^{-2 L_n} (1 - q^n( z_1, z_2 ) ) \\  + \frac{1}{2} e^{-2 L_n} \frac{1}{(V_{r_\alpha})^2} \int_{B(x,r_\alpha)^2} \left( q^n( y_1, y_2) - q^n(z_1, z_2 ) \right) dy_1 dy_2 - (1-e^{-2L_n}) q^n( z_1, z_2 ) \Bigg\rbrace dx = 0.
\end{multline*}
The first sum accounts for the contribution of events in which only one lineage is involved (hence the term $u^n_\alpha \left[ 1 - u^n_\alpha \1{B(z_j, r_\alpha)}(x) \right]$), and the last sum is the contribution of events in which both lineages are involved, which can result either in the coalescence of the two lineages, a simultaneous jump of the lineages, or in a recombination of at least one of the lineages in $]\ell, \ell + L_n[$.
We rewrite this equation in the more compact form
\begin{equation} \label{eq:qn}
    \mathcal{L}_2^n q^n ( z_1, z_2) - r_n(z_1, z_2) q^n( z_1, z_2) + g_n(z_1, z_2) = 0,
\end{equation}
where
\begin{equation*}
    g_n(z_1, z_2) := \frac{1}{2 N_n} e^{-2L_n} \sum_{\alpha \in \PM} \upsilon_\alpha \frac{\abs{\mathbb{H}_\alpha \cap B(z_1, r_\alpha) \cap B(z_2, r_\alpha)}}{V_{r_\alpha}},
\end{equation*}
$r_n$ is given by
\begin{multline*}
    r_n(z_1, z_2) := (1-e^{-L_n}) \sum_{i = 1}^2 \sum_{\alpha \in \PM} \frac{\abs{\mathbb{H}_\alpha \cap B(z_i, r_\alpha)}}{V_{r_\alpha}} \\ + \frac{1}{N_n} \left( \frac{1}{2} e^{-2 L_n} + (1 - e^{-2L_n}) - 2(1-e^{-L_n}) \right) \sum_{\alpha \in \PM} \upsilon_\alpha \frac{\abs{\mathbb{H}_\alpha \cap B(z_1, r_\alpha) \cap B(z_2, r_\alpha)}}{V_{r_\alpha}},
\end{multline*}
and $\mathcal{L}^n_2$ is the generator of a Markov process on $\R^2 \times \R^2$ given by
\begin{multline*}
    \mathcal{L}^n_2 f( z_1, z_2 ) := e^{-L_n}  \int_{\R^2} \left( \Phi(z_1, y) - \frac{1}{N_n} \Psi(z_1, z_2, y) \right) \left( f (y, z_2) - f( z_1, z_2 ) \right) dy \\ + e^{-L_n}  \int_{\R^2} \left( \Phi(z_2, y) - \frac{1}{N_n} \Psi(z_1, z_2, y) \right) \left( f (z_1, y) - f( z_1, z_2 ) \right) dy \\ + \frac{1}{2 N_n} e^{-2L_n} \int_{\R^2 \times \R^2} \Psi(z_1, z_2, y_1, y_2) \left( f( y_1, y_2 ) - f( z_1, z_2 ) \right) dy_1 dy_2,
\end{multline*}
where $\Phi$ was defined in \eqref{def:Phi},
\begin{equation*}
    \Psi(z_1, z_2, y) := \sum_{\alpha \in \PM} \upsilon_\alpha \frac{\abs{\mathbb{H}_\alpha \cap B(z_1, r_\alpha) \cap B(z_2, r_\alpha) \cap B(y, r_\alpha)}}{(V_{r_\alpha})^2},
\end{equation*}
and
\begin{equation*}
    \Psi(z_1, z_2, y_1, y_2) := \sum_{\alpha \in \PM} \upsilon_\alpha \frac{\abs{\mathbb{H}_\alpha \cap B(z_1, r_\alpha) \cap B(z_2, r_\alpha) \cap B(y_1, r_\alpha) \cap B(y_2, r_\alpha)}}{(V_{r_\alpha})^3}.
\end{equation*}

We can notice several things at once.
First, equation \eqref{eq:qn} is linear, with a source term $g_n$ which is of the order of $\frac{1}{N_n}$, which explains the fact that we need to multiply $p^n_{L_n, d_n}(\phi, \psi)$ by $N_n$ in \eqref{cvg_pn}.
In fact, the term $g_n$ corresponds to the rate of coalescence (with no recombination on $[\ell, \ell + L_n[$) of two lineages at locations $z_1$ and $z_2$.
The term $r_n$, on the other hand, is equal to $2 (1-e^{-L_n}) \sim 2 L_n$ as soon as the first coordinates of both $z_1$ and $z_2$ are larger than $R := \max(r_-, r_+)$, and when the distance between $z_1$ and $z_2$ is larger than $2R$.
Finally, the operator $\mathcal{L}^n_2$ encodes the reproduction events in which neither lineage recombines in $[\ell, \ell + L_n[$ and both lineages do not coalesce.
We see from the definition of $\mathcal{L}^n_2$ that, as long as the two lineages are further apart than $2R$, they behave like independent random walks with generator given by \eqref{def:generator_single_lineage} up to a time change of $e^{-L_n} \sim 1$.

Let $((\xi^{n,1}_t, \xi^{n,2}_t), t \geq 0)$ be a jump Markov process in $\R^2 \times \R^2$ with generator $\mathcal{L}^n_2$.
Then, by the Feynman-Kac formula, the solution to \eqref{eq:qn} can be expressed as
\begin{equation*}
    q^n(z_1, z_2) = \int_0^\infty \E[z_1, z_2]{ \exp\left( - \int_0^t r_n(\xi^{n,1}_s, \xi^{n,2}_s) ds \right) g_n(\xi^{n,1}_t, \xi^{n,2}_t) } dt.
\end{equation*}
Setting
\begin{equation*}
    X^{n,i}_t := \frac{1}{\sqrt{n}} \xi^{n,i}_{n t},
\end{equation*}
and
\begin{align*}
     \bar{r}_n (x_1, x_2) := n r_n(\sqrt{n} x_1, \sqrt{n} x_2), &&
     \bar{g}_n(x_1, x_2) := n g_n(\sqrt{n} x_1, \sqrt{n} x_2),
\end{align*}
we obtain by a change of variables in the integral over $t$,
\begin{equation*}
    q^n(\sqrt{n} x_1, \sqrt{n} x_2) = \int_0^\infty \E[x_1, x_2]{\exp\left( - \int_0^t \bar{r}_n(X^{n,1}_s, X^{n,2}_s) ds \right) \bar{g}_n(X^{n,1}_t, X^{n,2}_t) } dt.
\end{equation*}
Since the functions $\Phi$ and $\Psi$ are symmetric, the process $((X^{n,1}_t, X^{n,2}_t), t \geq 0)$ is reversible with respect to Lebesgue measure on $\R^2 \times \R^2$, and as a result,
\begin{multline} \label{reversibility}
    \int_{\R^2 \times \R^2} q^n(\sqrt{n} x_1, \sqrt{n} x_2) \phi(x_1) \psi(x_2) dx_1 dx_2 = \int_{\R^2 \times \R^2} \int_0^\infty \Bigg\lbrace \bar{g}_n(x_1, x_2)\\ \times \E[x_1, x_2]{\exp\left( - \int_0^t \bar{r}_n(X^{n,1}_s, X^{n,2}_s) ds \right) \phi(X^{n,1}_t) \psi(X^{n,2}_t) } \Bigg\rbrace dt dx_1 dx_2,
\end{multline}
which is also $p^n_{L_n, d_n}(\phi, \psi)$ by \eqref{eq_pn_qn}.

Recall that $X^{n,1}$ and $X^{n,2}$ behave approximately as independent copies of the rescaled process $X^n$ defined in Lemma~\ref{lemma:cvg_single_lineage}, which states that $X^n$ converges in distribution as $n \to \infty$ to the skew diffusion defined in \eqref{eds:skew_Bm}.
In addition to the fact that $\bar{r}_n$ is approximately $2 n L_n$ which converges to $2 \gamma$, we obtain the following.

\begin{lemma} \label{lemma:cvg_lineages}
Let $(\mu_n, n \geq 1)$ be a sequence of probability measures on $\R^2 \times \R^2$, which converges weakly as $n \to \infty$ to some probability measure $\mu$.
Let $\E[\mu_n]{\cdot}$ denote the expectation with respect to the law of $((X^{n,1}_t, X^{n,2}_t), t \geq 0)$, started from an initial configuration distributed according to $\mu_n$.
Then, under the above assumptions, as $n \to \infty$, for all $t \geq 0$,
\begin{multline*}
    \E[\mu_n]{\exp\left( - \int_0^t \bar{r}_n(X^{n,1}_s, X^{n,2}_s) ds \right) \phi(X^{n,1}_t) \psi(X^{n,2}_t) }  \\ \to e^{-2\gamma t} \int_{\R^2 \times \R^2} P_t\phi(x_1) P_t\psi(x_2) \mu(dx_1, dx_2).
\end{multline*}
\end{lemma}

Before proving this lemma, let us conclude the proof of \eqref{cvg_pn}.
To do so, we set
\begin{equation*}
    Q^n_{\phi, \psi}(x_1, x_2) := \int_0^\infty \E[x_1, x_2]{\exp\left( - \int_0^t \bar{r}_n(X^{n,1}_s, X^{n,2}_s) ds \right) \phi(X^{n,1}_t) \psi(X^{n,2}_t) } dt.
\end{equation*}
We then have, by \eqref{reversibility},
\begin{equation*}
    p^n_{L_n, d_n}(\phi, \psi) = \int_{\R^2 \times \R^2} \bar{g}_n(x_1, x_2) Q^n_{\phi, \psi}(x_1, x_2) dx_1 dx_2.
\end{equation*}
Replacing $\bar{g}_n$ by its definition yields
\begin{equation*}
    \frac{n}{2 N_n} e^{-2 L_n} \sum_{\alpha \in \PM} \frac{\upsilon_\alpha}{V_{r_\alpha}} \int_{\mathbb{H}_\alpha \times \R^2 \times \R^2} \1{\abs{x - \sqrt{n} x_1} < r_\alpha} \1{\abs{x - \sqrt{n} x_2} < r_\alpha} Q^n_{\phi, \psi}(x_1, x_2) dx_1 dx_2 dx.
\end{equation*}
Multiplying $x$ by $\sqrt{n}$ in the integral yields
\begin{equation*}
    \frac{e^{-2L_n}}{2 N_n} \sum_{\alpha \in \PM} \upsilon_\alpha V_{r_\alpha} \int_{\mathbb{H}_\alpha} \frac{1}{(V_{r_\alpha / \sqrt{n}})^2} \int_{B(x, r_\alpha / \sqrt{n})^2} Q^n_{\phi, \psi}(x_1, x_2) dx_1 dx_2 dx.
\end{equation*}
By Lemma~\ref{lemma:cvg_lineages}, as $n \to \infty$,
\begin{equation*}
    \frac{1}{(V_{r_\alpha / \sqrt{n}})^2} \int_{B(x, r_\alpha / \sqrt{n})^2} Q^n_{\phi, \psi}(x_1, x_2) dx_1 dx_2 \to \int_0^\infty e^{-2\gamma t} P_t\phi(x) P_t\psi(x) dt,
\end{equation*}
for any $x \in \R^2$.
As a result, recalling the definition of $\mathcal{N}(x)$ in \eqref{def:N}, and since $L_n \to 0$,
\begin{equation*}
    N_n p^n_{L_n, d_n}(\phi, \psi) \to \int_{\R^2} \int_0^\infty \frac{1}{\mathcal{N}(x)} e^{-2\gamma t} P_t\phi(x) P_t\psi(x) dt dx,
\end{equation*}
which proves \eqref{cvg_pn}.

\begin{proof}[Proof of Lemma~\ref{lemma:cvg_lineages}]
    Since $L_n \to 0$ as $n \to \infty$, we can safely ignore the factors $e^{-L_n}$ appearing in the definition of $\mathcal{L}_2$ to prove the result.
    We then note that $\mathcal{L}_2$ can be written (omitting the $e^{-L_n}$ factors) as
    \begin{equation}
        \mathcal{L}_2^n f(z_1, z_2) = \mathcal{L} f(\cdot, z_2)(z_1) + \mathcal{L} f(z_1, \cdot)(z_2) + \mathcal{R}^nf(z_1, z_2),
    \end{equation}
    where $\mathcal{L}$ is defined in \eqref{def:generator_single_lineage} and
    \begin{multline*}
        \mathcal{R}^n f(z_1, z_2) := \frac{1}{2 N_n} \int_{\R^2 \times \R^2} \Psi(z_1, z_2, y_1, y_2) (f(y_1, y_2) - f(z_1, z_2)) dy_1 dy_2 \\ - \frac{1}{N_n} \int_{\R^2} \Psi(z_1, z_2, y) (f(y, z_2) + f(z_1, y) - 2 f(z_1, z_2)) dy.
    \end{multline*}
    Also note that
    \begin{equation*}
        \int_{\R^2 \times \R^2} \Psi(z_1, z_2, y_1, y_2) dy_1 dy_2 = \int_{\R^2} \Psi(z_1, z_2, y) dy = \Psi(z_1, z_2).
    \end{equation*}
    We deduce that $(X^{n,1}, X^{n,2})$ can be coupled with a pair of independent copies of $X^n$ (defined in Lemma~\ref{lemma:cvg_single_lineage}), that we denote by $(\Tilde{X}^{n,1}, \Tilde{X}^{n,2})$, in such a way that, for some random time $\tau_c^n$, $X^{n,i}_t = \Tilde{X}^{n,i}_t$ for all $t \in [0,\tau_c^n[$, and that $\tau_c^n$ can be chosen such that
    \begin{equation*}
        \P{\tau_c^n > t} \geq \E{\exp\left( - \frac{C n}{N_n} \int_0^t \Psi(\sqrt{n} \Tilde{X}^{n,1}_s, \sqrt{n} \tilde{X}^{n,2}_s ) ds \right)},
    \end{equation*}
    for some constant $C > 0$ (we have accounted for the rescaling in the definition of $X^{n,i}$ by including the $n$ and $\sqrt{n}$ factors above).
    Note that there exists a constant $C > 0$ such that
    \begin{equation*}
        \Psi(\sqrt{n} z_1, \sqrt{n} z_2) \leq C \, \1{|z_1 - z_2| \leq \frac{2R}{\sqrt{n}}}.
    \end{equation*}
    In addition, looking at the definition of $\bar{r}_n$, we can note that there exists a constant $C > 0$ such that
    \begin{equation*}
        \abs{ n\bar{r}_n(\sqrt{n} z_1, \sqrt{n} z_2) - 2 \gamma } \leq C \left( \sum_{i=1}^2 \1{| z_i^{(1)} | \leq \frac{R}{\sqrt{n}}} + \frac{n}{N_n} \1{|z_1 - z_2| \leq \frac{2R}{\sqrt{n}}} \right) + \littleO{1}.
    \end{equation*}
    Hence, to prove Lemma~\ref{lemma:cvg_lineages}, it is sufficient to show that, as $n \to \infty$,
    \begin{equation} \label{eq:estimates_occupation_times}
        \int_0^t \1{|(\tilde{X}^{n,i}_s)^{(1)}| \leq \frac{R}{\sqrt{n}}} ds \to 0, \quad \text{and} \quad \frac{n}{N_n} \int_0^t \1{\abs{\tilde{X}^{n,1}_s - \tilde{X}^{n,2}_s} \leq \frac{2R}{\sqrt{n}}} ds \to 0,
    \end{equation}
    in probability for any $t > 0$.
    But, by a change of variables in the time integral,
    \begin{equation*}
        \int_0^t \1{|(X^{n}_s)^{(1)}| \leq \frac{R}{\sqrt{n}}} ds = \frac{1}{n} \int_0^{nt} \1{|\xi^{(1)}_s| \leq R} ds.
    \end{equation*}
    By Lemma~\ref{lem:occupation_boundary}, the expectation of the right hand side is $\bigO{\sqrt{\frac{t}{n}}}$, which yields the first convergence in \eqref{eq:estimates_occupation_times}.
    Similarly,
    \begin{equation*}
        \frac{n}{N_n} \int_0^t \1{\abs{\tilde{X}^{n,1}_s - \tilde{X}^{n,2}_s} \leq \frac{2R}{\sqrt{n}}} ds = \frac{1}{N_n} \int_0^{nt} \1{|\xi^1_s - \xi^2_s| \leq 2 R} ds,
    \end{equation*}
    where $\xi^1$ and $\xi^2$ are two independent copies of $\xi$ (defined in Lemma~\ref{lemma:cvg_single_lineage}).
    We then note that we can cover $\R^2$ by a countable set of balls of radius $\tilde{R}$ (which depends only on $R$) such that if $|\xi^1_s - \xi^2_s| \leq 2R$, then at least one of these balls contains both $\xi^1_s$ and $\xi^2_s$.
    Moreover this covering can be chosen such that any point in $\R^2$ is covered by at most $k$ balls, for some $k \in \N$.
    Let $(x_i, i \in I)$ denote the centre of balls of such a covering.
    Then
    \begin{align*}
        \E{ \1{|\xi^1_s - \xi^2_s| \leq 2R}} &\leq \sum_{i \in I} \P{ \xi_s \in B(x_i, \tilde{R}) }^2 \\
        &\leq k \sup_{x \in \R^2} \P{ \xi_s \in B(x, \tilde{R}) }.
    \end{align*}
    By \eqref{local_clt_bound}, the right hand side is bounded by $\frac{C}{s}$ for some constant $C > 0$.
    It follows that
    \begin{equation*}
        \E{ \int_0^{nt} \1{|\xi^1_s - \xi^2_s| \leq 2R} ds } \leq C \log(1+nt),
    \end{equation*}
    for some constant $C > 0$.
    Hence the second part of \eqref{eq:estimates_occupation_times} follows by \eqref{parameter_scaling}.
    This concludes the proof of the lemma.
\end{proof}

\section{Convergence of ancestral lineages to skew Brownian motion} \label{app:cvg_to_skew_Bm}

The aim of this section is to prove Lemma~\ref{lemma:cvg_single_lineage}.
In \cite{iksanov_functional_2016}, the authors prove a similar result in the case of a discrete time Markov chain taking values in $\Z$ which behaves like a simple random walk outside a bounded region around the origin.
We thus extend their proof to continuous time jump Markov processes with continuous state space.
To simplify notations, we restrict the proof to the convergence of the first coordinate of $X^n$, but the generalisation to the two-dimensional case is straightforward.
In the following, we thus use $\xi_t$, $X^n_t$, etc. instead of $\xi^{(1)}_t$, $(X^n_t)^{(1)}$, etc.

Recall that we have set $R = \max(r_-, r_+)$.
For $\alpha \in \PM$, set
\begin{equation*}
\xi_\alpha (t) := \alpha \xi_t \, \1{ \alpha \xi_t > R },
\end{equation*}
and
\begin{align*}
\tau^\alpha_0 &:= \inf \lbrace t > 0 : |\xi_t| \leq R \rbrace, \\
\eta^\alpha_k &:= \inf \left\lbrace t > \tau^\alpha_k : \alpha \xi_t > R \right\rbrace, \quad k\geq 0, \\
\tau^\alpha_{k+1} &:= \inf \left\lbrace t > \eta^\alpha_k : \alpha \xi_t \leq R \right\rbrace, \quad k\geq 0.
\end{align*}
One can then write the decomposition (see formula (2.1) in \cite{iksanov_functional_2016})
\begin{equation} \label{decomposition}
\xi_\alpha(t) = \xi_\alpha(0) + M_\alpha(t) + L_\alpha(t) - \alpha \sum_{i \geq 0} \xi_{\tau_i^\alpha} \1{ \tau_i^\alpha \leq t < \eta_i^\alpha },
\end{equation}
with
\begin{align*}
M_\alpha(t) &:= \alpha \int_{0}^{t} \1{ \alpha \xi_{s^-} > R } d\xi_s, \\
L_\alpha(t) &:= \alpha \sum_{i\geq 0} \left( \xi_{\eta_i^\alpha} - \xi_{\tau_i^\alpha} \right) \1{ \eta_i^\alpha \leq t }.
\end{align*}
Here, $ M_\alpha(t) $ is the sum of the jumps of $ \xi $ from $ \{ x \in \R : \alpha x > R \} $ up to time $ t $.
Note that these jumps are all independent centred random variables, so $ (M_\alpha(t), t \geq 0) $ is a martingale with respect to the natural filtration of $ (\xi_t, t \geq 0) $.
On the other hand, $ L_\alpha(t) $ is an increasing process which only increases when $ \xi $ escapes from $ \{ \alpha x \leq R \} $, and should be thought of as an analogue of the left (for $ L_- $) or right (for $ L_+ $) local time of $ \xi $ at zero.
Also set
\begin{align*}
M^n_\alpha(t) := \frac{1}{\sqrt{n}} M_\pm( nt ), && L^n_\alpha(t) := \frac{1}{\sqrt{n}} L_\alpha( nt ).
\end{align*}
We also set, for $\alpha \in \PM$, $X^{n,\alpha}_t := \max(\alpha X^n_t, 0)$.
The following then holds.

\begin{lemma} \label{lem:tightness}
For any fixed $T>0$, the sequence of $\sko{\R^6}$-valued random variables 
\begin{equation} \label{eq:tight_sequence}
    \lbrace (X^{n,\alpha}, M_\alpha^n, L_\alpha^n), \alpha \in \PM, n \in \N \rbrace
\end{equation}
is tight.
Furthermore, any limit point $ \lbrace ( X^{\infty}_{\alpha}, M_\alpha^\infty, L_\alpha^\infty ), \alpha \in \PM \rbrace$ of the sequence is a continuous process satisfying
\begin{equation} \label{zero_occupation_time}
\int_{0}^{T} \1{ X^{\infty}_+(t) - X^{\infty}_-(t) = 0} dt = 0, \quad \text{a.s.}
\end{equation}
\end{lemma}

\begin{lemma} \label{lem:qvar}
Let $\lbrace ( X^{\infty}_{\alpha}, M_\alpha^\infty, L_\alpha^\infty ), \alpha \in \PM \rbrace$ be the limit point of a converging subsequence of \eqref{eq:tight_sequence}.
Then
\begin{enumerate}[i)]
\item the processes $L^\infty_\alpha$ are non-decreasing almost surely and satisfy
\begin{equation*}
\int_{0}^{T} \1{ X^{\infty}_{\alpha}(t) > 0} d L_\alpha^\infty (t) = 0 \quad \text{a.s.}
\end{equation*}
\item the processes $M_\alpha^\infty$ are continuous martingales with respect to the filtration generated by $(X^{\infty}_{\alpha}(t), L_\alpha^\infty(t), M_\alpha^\infty(t), t \geq 0)$ with predictable quadratic variation
\begin{equation*}
\langle M_\alpha^\infty \rangle_t = \sigma_\alpha^2 \int_{0}^{t} \1{ X^{\infty}_{\alpha} (s) > 0} ds.
\end{equation*}
\end{enumerate}
\end{lemma}

\begin{lemma} \label{lem:local_times}
Let $\lbrace ( X^{\infty}_{\alpha}, M_\alpha^\infty, L_\alpha^\infty ), \alpha \in \PM \rbrace$ be the limit point of a converging subsequence of \eqref{eq:tight_sequence}.
Then for any $t\geq 0$, 
\begin{equation*}
L_+^\infty(t) = \frac{\sigma_+^2}{\sigma_-^2} L_-^\infty(t),
\end{equation*}
almost surely.
\end{lemma}

Lemma~\ref{lemma:cvg_single_lineage} follows from the above lemmas and Proposition~2.1 in \cite{iksanov_functional_2016}.
Lemma~\ref{lem:tightness} is proved in Subsection~\ref{subsec:tightness}.
The proof of Lemma~\ref{lem:qvar} does not differ significantly from the one given for Lemma~2.2 in \cite{iksanov_functional_2016} and we omit the details.
The proof of Lemma~\ref{lem:local_times} is given in Subsection~\ref{subsec:local_time}.

\subsection{Occupation time of the boundary} \label{subsec:occupation_time}

We begin with the following result controlling the time spent by $\proc{\xi}$ in the region $[-R,R]$.
This lemma is used in the proof of Lemma~\ref{lem:qvar} (see the proof of Lemma~2.3 in \cite{iksanov_functional_2016}) and in the proof of Lemma~\ref{lem:local_times}.

\begin{lemma} \label{lem:occupation_boundary}
	For $t \geq 0$, define $\nu(t) = \int_{0}^{t} \1{\abs{\xi_s} \leq R} ds $.
	Then
	\begin{enumerate}[i)]
		\item $\lim_{t\to\infty} \nu(t) = +\infty$ almost surely,
		\item $\sup_{x\in\R} \E[x]{ \nu(t) } = \bigO{\sqrt{t}}$ a.s. as $t\to\infty$.
	\end{enumerate}
\end{lemma}

\begin{proof}
	Since $ (\xi_t)_{t \geq 0} $ is neighbourhood-recurrent, $ \nu(t) \to \infty $ as $ t \to \infty $.
	Set $\zeta_0 = 0$ and
	\begin{align*}
	& \varsigma_i := \inf \left\lbrace t > \zeta_{i-1} : \abs{\xi_t} \leq R \right\rbrace, \quad i \geq 1, \\
	& \zeta_i := \inf \left\lbrace t > \varsigma_i : \abs{\xi_t} > R \right\rbrace, \quad i \geq 1.
	\end{align*}
	Then $\nu(t)$ can be written as the sum of the lengths of the excursions inside $[-R,R]$ up to time $t$,
	\begin{equation*}
	\nu(t) = \sum_{i\geq 1} \left( \zeta_i \wedge t - \varsigma_i \wedge t \right).
	\end{equation*}
	Hence
	\begin{equation*}
	\E[x]{ \nu(t) } \leq \E[x]{ \sum_{i\geq 1} \E{ \zeta_i - \varsigma_i }{ \F_{\varsigma_i} } \1{\varsigma_i \leq t} }.
	\end{equation*}
	Noting that, since
	\begin{equation*}
	    \int_{|y| > R} \Phi(x,y) dy > 0, \quad \forall\, x \in [-R,R],
	\end{equation*}
	there exists $\varepsilon>0$ such that, for all $ dt > 0 $ small enough, 
	\begin{equation*}
	    \P{ \abs{\xi(t + dt)} > R }{ \xi_t=x } = \frac{1}{\lambda(x)} \geq \varepsilon dt
	\end{equation*}
	for all $\abs{x} \leq R$, we see that $\zeta_i - \varsigma_i$ is stochastically dominated by an exponential random variable with parameter $\varepsilon$.
	Hence
	\begin{equation*}
	\E[x]{ \nu(t) } \leq \frac{1}{\varepsilon} \E[x]{ \sum_{i \geq 1} \1{\varsigma_i \leq t} }.
	\end{equation*}
	In addition, the number of visits to $[-R,R]$ before time $t$ is less than the number of visits to this set before the first excursion out of this set longer than $t$, \textit{i.e.}
	\begin{equation*}
	\sum_{i \geq 1} \1{\varsigma_i \leq t} \leq m(t) := \inf \lbrace i \geq 1 : \varsigma_{i+1} - \zeta_i > t \rbrace.
	\end{equation*}
	Let $\proc{W}$ be a continuous time random walk on $ \R $ with infinitesimal generator
	\begin{align*}
		\mathcal{G} f(x) = \int_{\R} (f(y) - f(x)) \frac{\abs{B(x,1) \cap B(y,1)}}{V_{1}^2} dy.
	\end{align*}
	Then for any $ x > R $,
	\begin{equation*}
		\P[\alpha x]{\varsigma_1 - \zeta_0 > t} \geq \P[0]{\inf_{0 \leq s \leq t} W_s \geq 0}.
	\end{equation*}
	(Notice that the right hand side is unchanged if $W$ is replaced by $R W$.)
	As a result $m(t)$ is stochastically dominated by a geometric random variable with parameter
	\begin{equation*}
	p(t) = \P[0]{ \inf_{0 \leq s \leq t} W_s \geq 0 }.
	\end{equation*}
	Furthermore, there exists $\eta > 0$ such that, for all $t \geq 0$, $ p(t) \geq \frac{\eta}{\sqrt{t}} $, (see pp. 381-382 in \cite{bingham_regular_1989} or equations (3.4) and (3.5) in \cite{iksanov_functional_2016}).
	As a result, for all $ x \in \R $,
	\begin{equation*}
	\E[x]{ \nu(t) } \leq \frac{1}{\varepsilon p(t)} \leq \frac{\sqrt{t}}{\varepsilon \eta}.
	\end{equation*}
\end{proof}

\subsection{Tightness} \label{subsec:tightness}

Let us now give the proof of Lemma~\ref{lem:tightness}.

\begin{proof}[Proof of Lemma~\ref{lem:tightness}]
From \eqref{decomposition}, and the fact that
\begin{equation*}
	\Bigg| \sum_{i\geq 0} \xi_{\tau_i^\alpha} \1{\tau_i^\alpha \leq t < \eta_i^\alpha} \Bigg| \leq R,
\end{equation*}
it is enough to prove the tightness of $X^n$ and $M^n_\alpha$.
We use Aldous' criterion \cite{aldous_stopping_1978} to prove that $ M^n_\alpha $ is tight and then we use the fact that the increments of $ \xi $ are bounded by those of
\begin{equation*}
    M := M_+ - M_-
\end{equation*}
(see equation \eqref{bound_increments} below) to show that $ X^n $ (and hence $X^{n,\alpha}$) is tight.

From the definition of $\xi$, we have $M^n_\alpha(0) = 0$ and
\begin{equation*}
\sup_{t\geq 0} \abs{ M^n_\alpha(t) - M^n_\alpha(t_-) } \leq \frac{2R}{\sqrt{n}}.
\end{equation*}
Moreover, since outside $[-R,R]$, $\xi$ behaves as a simple random walk, for any stopping time $S$ and $\delta > 0$,
\begin{equation*}
\E{ \left( M^n_\alpha( S + \delta ) - M^n_\alpha ( S ) \right)^2 } \leq \sigma_\alpha^2 \delta,
\end{equation*}
proving that $(M^n_\alpha, n \in \N)$ is C-tight by \citep[Theorem~1]{aldous_stopping_1978}.

Now take $ 0 \leq s \leq t $. 
If $ \xi $ does not visit $ [-R,R] $ between time $ s $ and time $ t $, then $ \xi_t - \xi_s = M(t) - M(s) $.
If it does visit this set, then let $ t_1 $ be the first time $ \xi $ enters $ [-R,R] $ after time $ s $ and $ t_2 $ the last time $ \xi $ leaves this set before time $ t $.
Then
\begin{align*}
	\abs{ \xi_t - \xi_s } &\leq \abs{\xi_t - \xi_{t_2}} + \abs{\xi_{t_2} - \xi_{t_1}} + \abs{\xi_{t_1} - \xi_s} \\
	&\leq 2 R + \abs{M(t) - M(t_2)} + \abs{M(t_1) - M(s)}.
\end{align*}
As a result, for $\delta > 0$,
\begin{equation} \label{bound_increments}
\sup_{\abs{s-t} \leq \delta n} \abs{ \xi_s - \xi_t } \leq 2R + 2 \sup_{\abs{s-t} \leq \delta n} \abs{ M(s) - M(t) }.
\end{equation}
The tightness of $(X^n, n \in \N)$ then follows from that of $(M^n, n \in \N)$ by writing
\begin{multline} \label{tightness_X}
\lim_{\delta \downarrow 0} \limsup_{n \to \infty} \mathbb{P}\Big( \sup_{\substack{ \abs{s-t} \leq \delta n \\ s, t \in [0, n T] }} \abs{ \xi_s - \xi_t } > \varepsilon \sqrt{n} \Big) \\
\leq \lim_{\delta \downarrow 0} \limsup_{n \to \infty} \mathbb{P}\Big( 2R + 2 \sup_{\substack{ \abs{s-t} \leq \delta n \\ s, t \in [0, n T] }} \abs{ M(s) - M(t) } > \varepsilon \sqrt{n} \Big) = 0.
\end{multline}
It remains to prove \eqref{zero_occupation_time}.
Note that any limit point $ \lbrace ( X_\alpha^\infty, M_\alpha^\infty, L_\alpha^\infty ), \alpha \in \PM \rbrace $ satisfies
\begin{equation*}
X^\infty(t) = X^\infty_+(t) - X_-^\infty(t) = M_+^\infty(t) - M_-^\infty(t) + L_+^\infty(t) - L_-^\infty(t) = M^\infty(t) + L^\infty(t).
\end{equation*}
From the definition of $M^n_\alpha$ and Lemma~\ref{lem:occupation_boundary}, one shows, as in the proof of Lemma~2.3 in \cite{iksanov_functional_2016}, that $M^\infty$ is a stochastic integral with respect to standard Brownian motion $(B_t, t \geq 0)$,
\begin{equation*}
M^\infty(t) = \int_{0}^{t} \sigma( X^\infty(s) ) dB_s.
\end{equation*}
In addition, $L_\alpha^\infty$ is a continuous process with locally bounded variation.
As a result $\langle X^\infty \rangle_t = \langle M^\infty \rangle_t$ and \eqref{zero_occupation_time} follows from the occupation density formula.
\end{proof}

\subsection[The left and right local time at zero of xi]{The left and right local time at zero of $\proc{\xi}$} \label{subsec:local_time}

The proof of Lemma~\ref{lem:local_times} is adapted from that of Lemma~2.3 in \cite{iksanov_functional_2016}.
Recall the expression for the left and right local time of $ \proc{\xi} $,
\begin{equation*}
	L_\alpha(t) = \alpha \sum_{i \geq 0} (\xi_{\eta_i^\alpha} - \xi_{\tau_i^\alpha}) \1{\eta_i^\alpha \leq t}.
\end{equation*}
For any particular visit of $ \xi $ to $ [-R,R] $, the value of $ \xi_{\eta_i^\alpha} - \xi_{\tau_i^\alpha} $ depends on the value of $ \xi $ when it enters this set.
However, over many visits to $ [-R,R] $, $ L_\alpha(t) $ only records an average of these values.
The "typical" value of $ \xi_{\eta_i^\alpha} - \xi_{\tau_i^\alpha} $ can thus be expressed with the help of the stationary distribution of the process describing the visits of $ \xi $ to $ [-R,R] $ ($ (Y(t), t \geq 0) $ below).
The left and right local time of $ \xi $ then become asymptotically proportional to the occupation time of the boundary $ \nu(t) = \int_{0}^{t} \1{\abs{\xi_s} \leq R} ds $, with different coefficients whose expressions we give below.

Set, for $ t\geq 0 $,
\begin{equation*}
	\theta(t) = \inf \lbrace \theta > 0 : \nu(\theta) > t \rbrace, \quad Y(t) := \xi_{\theta(t)}.
\end{equation*}
The process $ (Y(t), t \geq 0) $ is a jump Markov process taking values in $ [-R,R] $, describing the values taken by $ \xi $ inside this region.
Let $ \bar{\theta} $ denote the left-continuous version of $ \theta $, \textit{i.e.} for $ t\geq 0 $,
\begin{equation*}
\bar{\theta}(t) = \sup \lbrace \theta \geq 0 : \nu(\theta) < t \rbrace.
\end{equation*}
If $ t\geq 0 $ is such that $ \bar{\theta}(t) \neq \theta(t) $, then $ \xi $ makes an excursion outside $ [-R,R] $ in the interval $ [\bar{\theta}(t), \theta(t)] $.

For $ \alpha \in \PM $, let $ V_\alpha $ be defined by
\begin{equation*}
V_\alpha (t) := \alpha \left( Y(t) - Y(0) \right) + \alpha \sum_{0 < s\leq t} (\xi_{\bar{\theta}(s)} - \xi_{\theta(s))} \1{ \alpha \xi_{\bar{\theta}(s)} > R }.
\end{equation*}
Then $ (V_\alpha(t), t \geq 0) $ is a process which has the two following nice properties.

\begin{lemma} \label{lemma:V_approx_L}
	For all $ t \geq 0 $,
	\begin{align*}
		\abs{V_\alpha(\nu(t)) - L_\alpha(t)} \leq 4 R.
	\end{align*}
\end{lemma}

\begin{lemma} \label{lemma:V_martingale}
	The process $ (Z_\alpha(t), t \geq 0) $ defined by
	\begin{align}\label{def_Zt}
		Z_\alpha(t) := V_\alpha(t) - \int_{0}^{t} h_\alpha(Y(s)) ds
	\end{align}
	with
	\begin{align} \label{function_h}
		h_\alpha(x) := \alpha \int_\R \Phi(x,y) \1{\alpha y \leq R} (\E[y]{\xi_{\theta(0)}} - x) dy + \alpha \int_\R \Phi(x,y) \1{\alpha y > R} (y-x) dy
	\end{align}
	is a square-integrable martingale with respect to the filtration $ (\F_{\theta(t)}, t \geq 0) $ (where $ (\F_t, t \geq 0) $ is the natural filtration of $ (\xi_t, t \geq 0) $).
	Furthermore, for all $ t \geq 0 $,
	\begin{align*}
		\langle Z_\alpha \rangle_t \leq 4 R^2 t.
	\end{align*}
\end{lemma}

These two lemmas are proved in Subsection~\ref{subsec:process_V}.
In addition, we have the following.

\begin{lemma} \label{lemma:stationary_dist_Y}
	The uniform probability measure on $[-R,R]$ is a stationary measure for the Markov process $ (Y(t), t \geq 0) $ and $Y$ is ergodic with respect to this measure.
	Moreover,
	\begin{align*}
		\frac{1}{2 R}\int_{[-R,R]} h_\alpha(x) dx = \frac{\sigma_\alpha^2}{4 R}.
	\end{align*}
\end{lemma}

Lemma~\ref{lemma:stationary_dist_Y} is proved in Section~\ref{sec:proof_lemma_h} of the supplemental material.
With these results, we now prove Lemma~\ref{lem:local_times}.

\begin{proof}[Proof of Lemma~\ref{lem:local_times}]
	First note that, by Lemma~\ref{lemma:V_martingale} and Theorem~1 in \cite{lepingle_sur_1978},
	\begin{align*}
		\frac{1}{t} V_\alpha(t) - \frac{1}{t} \int_{0}^{t} h_\alpha(Y(s)) ds \cvgas{t} 0
	\end{align*}
	almost surely.
	Then, by Lemma~\ref{lemma:stationary_dist_Y} and the pointwise ergodic theorem,
	\begin{align*}
		\frac{1}{t} \int_{0}^{t} h_\alpha(Y(s)) ds \cvgas{t} \frac{1}{2R} \int_{[-R,R]} h_\alpha(x) dx,
	\end{align*}
	almost surely.
	From Lemma~\ref{lemma:V_approx_L} and Lemma~\ref{lem:occupation_boundary}.i, we obtain that
	\begin{align*}
		\frac{1}{\nu(t)} L_\alpha(t) \cvgas{t} \frac{1}{2R} \int_{[-R,R]} h_\alpha(x) dx,
	\end{align*}
	almost surely, and Lemma~\ref{lemma:stationary_dist_Y} then implies
	\begin{align*}
		\lim_{t \to \infty} \frac{L_+(t)}{L_-(t)} = \frac{\sigma_+^2}{\sigma_-^2}.
	\end{align*}
\end{proof}

\subsection[The process Valpha(t), t > 0]{The process $ (V_\alpha(t), t \geq 0) $} \label{subsec:process_V}

\begin{proof}[Proof of Lemma~\ref{lemma:V_approx_L}]
	Note that $ \alpha \xi_{\bar{\theta}(s)} > R $ with $ s>0 $ if and only if $ \bar{\theta}(s) = \eta^\alpha_i $ for some $ i \geq 0 $, and in this case, $ \theta(s) = \tau^\alpha_{i+1} $.
	In addition, $ s \leq \nu(t) $ if and only if $ \bar{\theta}(s) \leq t $, as a result,
	\begin{equation*}
	V_\alpha (\nu(t)) = \alpha (Y(\nu(t)) - Y(0)) + \alpha \sum_{i\geq 0} (\xi_{\eta^\alpha_i} - \xi_{\tau^\alpha_{i+1}}) \1{\eta^\alpha_i \leq t}.
	\end{equation*}
	Hence
	\begin{align*}
	\abs{V_\alpha(\nu(t)) - L_\alpha(t)} \leq \abs{Y(\nu(t))} + \abs{Y(0)} + \abs{\xi_{\tau^\alpha_0}} + \sum_{i \geq 2} \abs{\xi_{\tau^\alpha_i}} \1{\eta^\alpha_{i-1} \leq t < \eta^\alpha_i}.
	\end{align*}
	Since $ \abs{Y(t)} \leq R $, $ \abs{\xi_{\tau_i^\alpha}} \leq R $,
	\begin{equation*}
	\abs{ V_\alpha(\nu(t)) - L_\alpha(t) } \leq 4 R.
	\end{equation*}
\end{proof}

\begin{proof}[Proof of Lemma~\ref{lemma:V_martingale}]
	We start by introducing the following notation.
	For $ f : [-R,R] \to \R $ measurable and bounded and $ x \in \R $, let
	\begin{equation} \label{definition_operator_E}
	Ef(x) := \E[x]{f(Y(0))} = \E[x]{f(\xi_{\theta(0)})}
	\end{equation}
	(recall that $ \theta(0) = \inf \{ t \geq 0 : \abs{\xi_t} \leq R \} $).
	Also let $ \iota : [-R,R] \to \R $ be defined by $ \iota(x) = x $.
	We then note that $ (Y(t), V_\alpha(t))_{t \geq 0} $ is a jump Markov process with respect to the filtration $ (\F_{\theta(t)}, t \geq 0) $ with generator
	\begin{multline*}
		\mathcal{G}_\alpha f(x,v) := \int_\R \Phi(x,y) \1{\alpha y \leq R} \left( f(E \iota(y), v + \alpha (E\iota(y) - x)) - f(x,v) \right) dy \\ + \int_\R \Phi(x,y) \1{\alpha y > R} \left( f(E\iota(y), v + \alpha (y-x)) - f(x,v) \right) dy.
	\end{multline*}
	Noting that, for $ f(x,v) = v $,
	\begin{align*}
		\mathcal{G}_\alpha f(x,v) = h_\alpha(x),
	\end{align*}
	we obtain that $ (Z_\alpha(t), t \geq 0) $ defined in \eqref{def_Zt} is a local martingale with respect to this filtration and that its predictable quadratic variation is
	\begin{align*}
		\langle Z_\alpha \rangle_t = \int_{0}^{t} Q_\alpha(Y(s)) ds
	\end{align*}
	with
	\begin{align*}
		Q_\alpha(x) := \int_\R \Phi(x,y) \1{\alpha y \leq R} (E\iota(y) - x)^2 dy + \int_\R \Phi(x,y) \1{\alpha y > R} (y-x)^2 dy.
	\end{align*}
	We then conclude by noting that $ \abs{Q_\alpha(x)} \leq 4 R^2 $.
\end{proof}

\section{Random walk approximation of the skew diffusion} \label{sec:random_walk}

\begin{figure}[hb]
    \centering
    \includegraphics[width=0.6\textwidth]{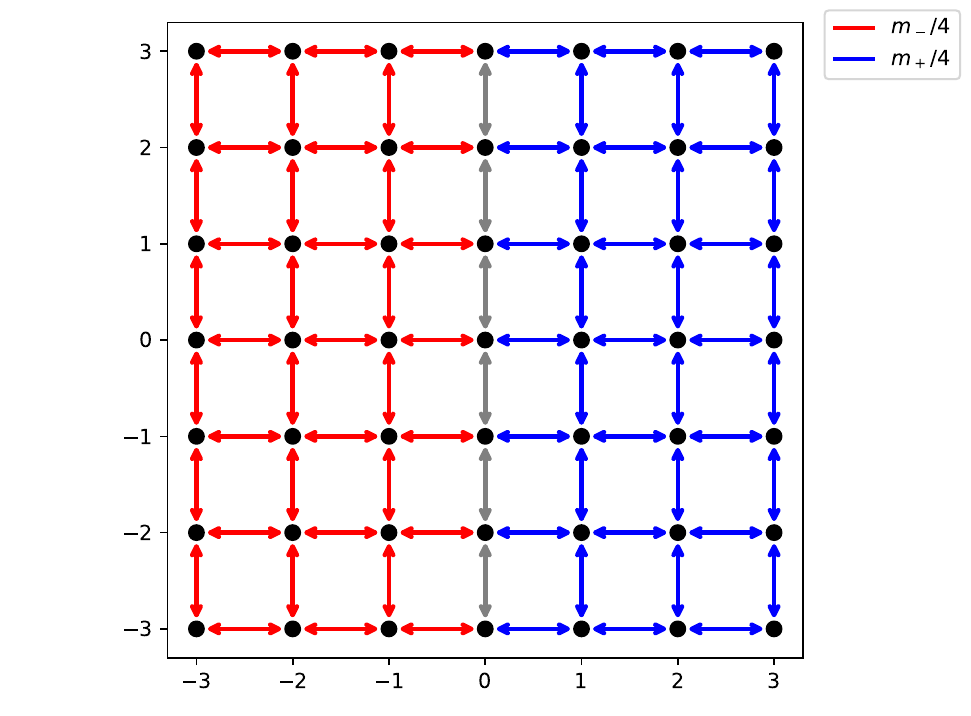}
    \caption{Transition kernel of the random walk $(Y^{(L)}_t, t \geq 0)$ used to approximate the transition density of the skew diffusion. At each step, the probability of jumping along each edge depends on the color of the edge ($m_-/4$ for red edges, $m_+/4$ for blue edges, and $(m_++m_-)/8$ for grey edges).}
    \label{fig:kernel}
\end{figure}

Fix $m_+ \in (0,1)$, $m_- \in (0,1)$ and $L \in 2 \N + 1$.
Consider a random walk $(Y_n^{(L)}, n \in \N)$ taking values in $\mathcal{X}_L := \left\llbracket -\frac{L + 1}{2}, \frac{L + 1}{2} \right\rrbracket^2$, with transition kernel given by
\begin{equation*}
    \Pi_{x,y} := \begin{cases}
                    \frac{m_\alpha}{4} &\text{ for } |x - y| = 1 \text{ and either } \alpha x^{(1)} > 0 \text{ or } \alpha y^{(1)} > 0, \\
                    \frac{m_+ + m_-}{8} &\text{ for } |x-y| = 1 \text{ and } x^{(1)} = y^{(1)} = 0, \\
                    1 - \sum_{z : |z - x| = 1} \Pi_{x, z} &\text{ for } y = x.
                \end{cases}
\end{equation*}
This kernel is depicted on Figure~\ref{fig:kernel}.
The following is then a direct consequence of \citep[Theorem~1.1]{iksanov_functional_2016}.

\begin{theorem} \label{thm:cvg_random_walk}
Suppose that $(L_n, n \in \N)$ and $(\delta_n, n \in \N)$ are two sequences such that $\delta_n \to 0$ and $ \delta_n L_n \to \infty$ as $n \to \infty$, and fix $r > 0$. 
Set
\begin{equation*}
    Y^n_t := \delta_n Y^{(L_n)}_{r t / \delta_n^2}, \quad t \geq 0.
\end{equation*}
Then, as $n \to \infty$, $(Y^n_t, t \geq 0)$ converges in distribution to $(X_t, t \geq 0)$, solution of \eqref{eds:skew_Bm}, with $\sigma_\alpha^2 = r \frac{m_\alpha}{2}$.
\end{theorem}

Given a set of sampling positions, we thus chose a step size $\delta$ and a grid size $L$ such that $\delta$ is small compared to the mean distance between samples and $L$ is large compared to the maximum distance between the samples.
We then choose $m_+$, $m_-$ and $r$ as functions of $\sigma_+$ and $\sigma_-$, and we approximate $P_t\phi(x)$ by
\begin{equation*}
    \E[x]{\phi(\delta \, Y^{(L)}_{t / \delta^2})}
\end{equation*}
in \eqref{eq:WM_E_N_geq_L}.
Note that if $\sigma_\alpha^2 < 1/2$ for $\alpha \in \PM$, then we can take $r = 1$, $m_\alpha = 2 \sigma_\alpha^2$ above, but if $\sigma_\alpha^2 \geq 1/2$ for some $\alpha$, then we choose $r > 0$ large enough that $\sigma_\alpha^2 / r < 1/2$ for $\alpha \in \PM$ and we take $m_\alpha = 2 \sigma_\alpha^2 / r$.

\begin{figure}[p]
\centering
\includegraphics[width=0.9\textwidth]{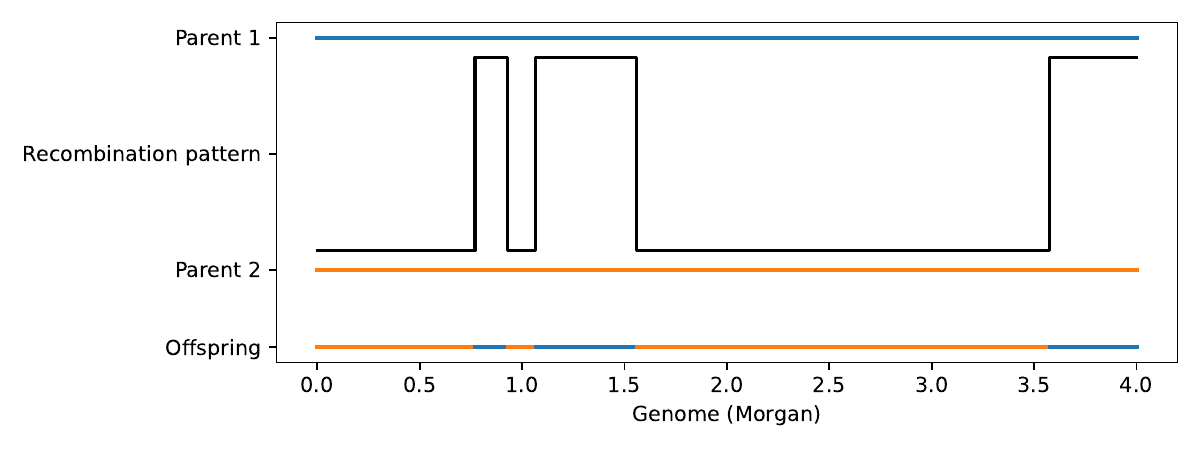}
\caption{Example of a recombination pattern and the resulting offspring genome. The crossover take place at the points of a standard Poisson process along the genome, and the parent genome that is ancestral to that of the new genome changes at each of these points.}\label{fig:recombination}
\end{figure}

\clearpage

\begin{figure}
    \centering
    \includegraphics[width=\textwidth]{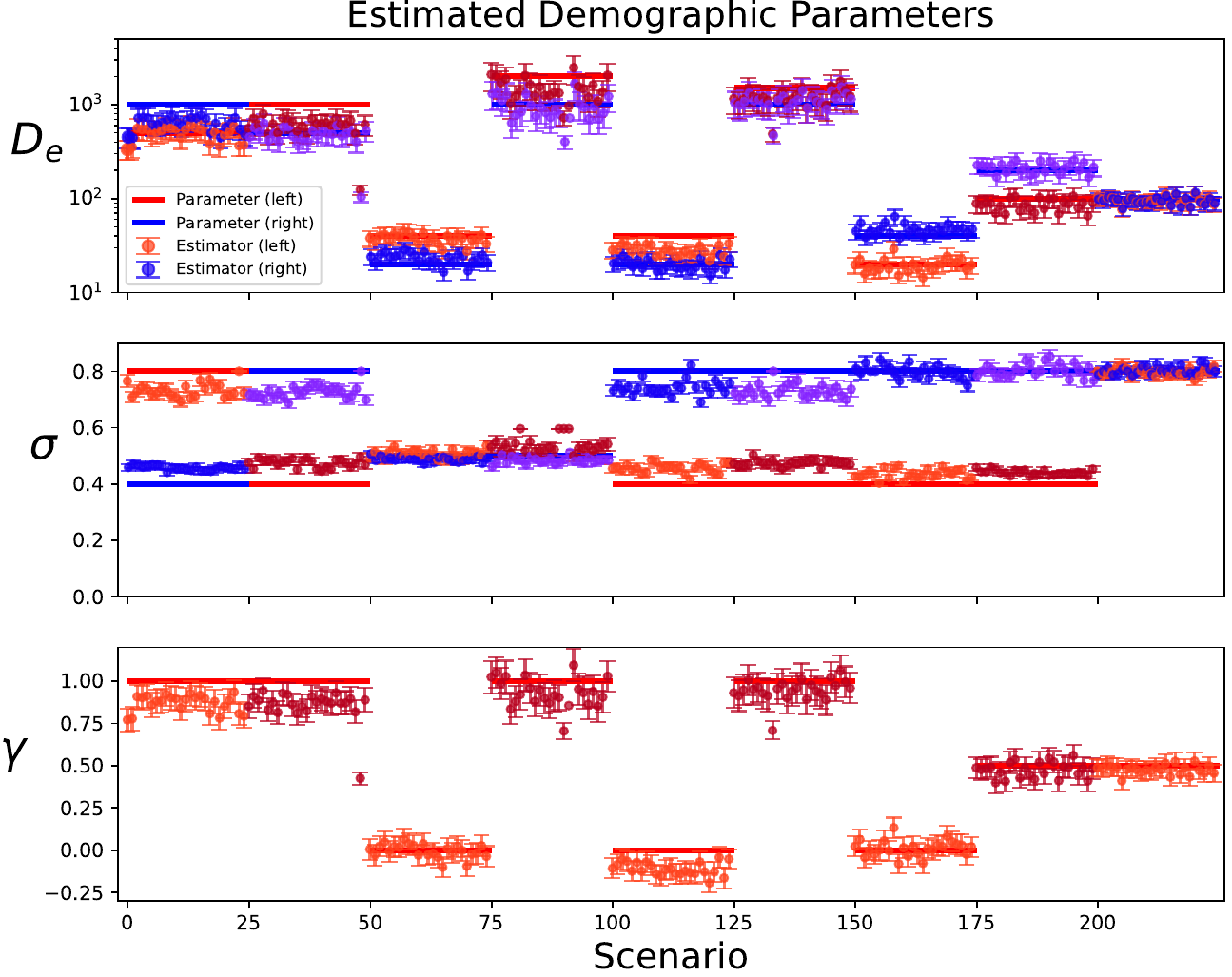}
    \caption{Simulated (straight lines) vs. estimated (dots) parameters for 20 independent runs for each of 9 scenarios. Diffusion and effective neighbourhood size on the left of the interface are shown in red, and those on the right are shown in blue. The error bars correspond to the 95\% confidence interval obtained from the Fisher information matrix (see the discussion in the main text). The parameter $\gamma$ is the growth exponent of the population density (it is the same on both sides of the interface), such that the effective neighbourhood size $t$ generations in the past is $D_e t^{-\gamma}$.}
    \label{fig:simu_results_varying_pop_size}
\end{figure}

\clearpage

\begin{figure}[p]
    \centering
    \includegraphics[width=0.9\textwidth]{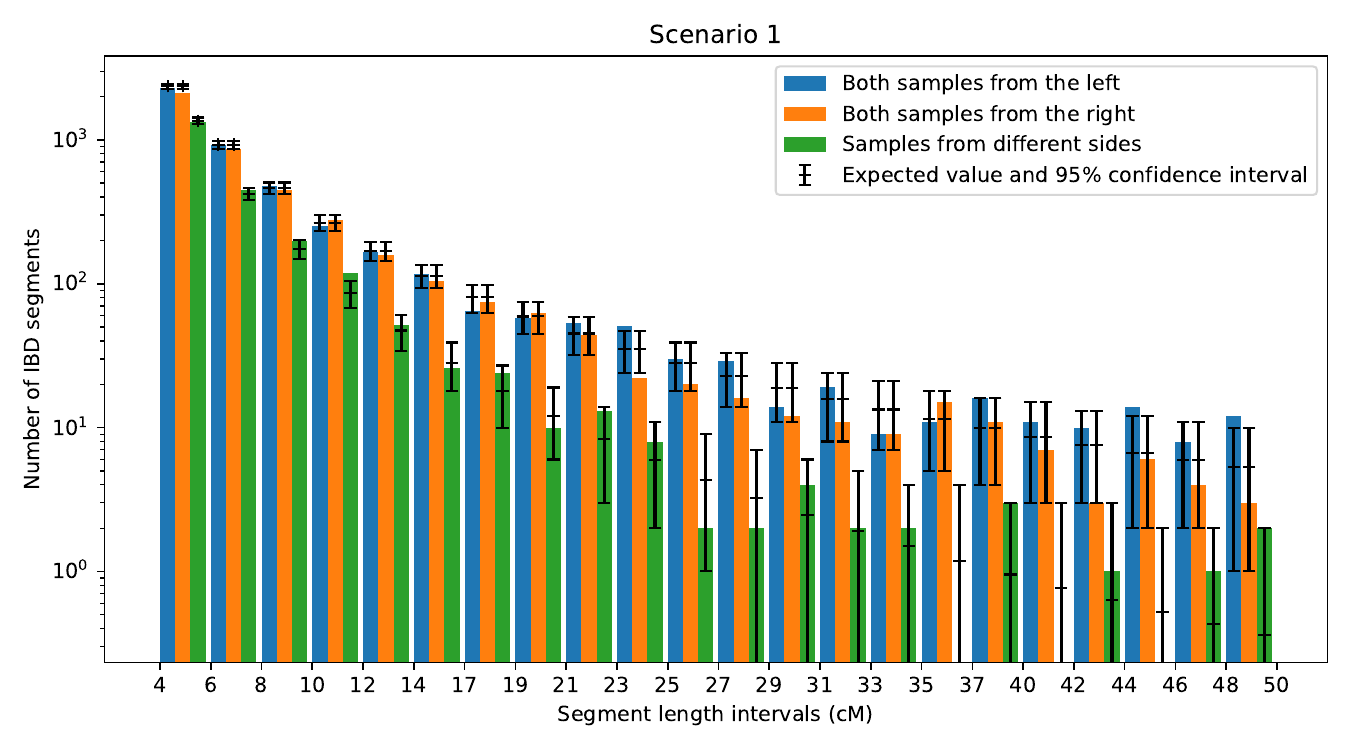}
    \includegraphics[width=0.9\textwidth]{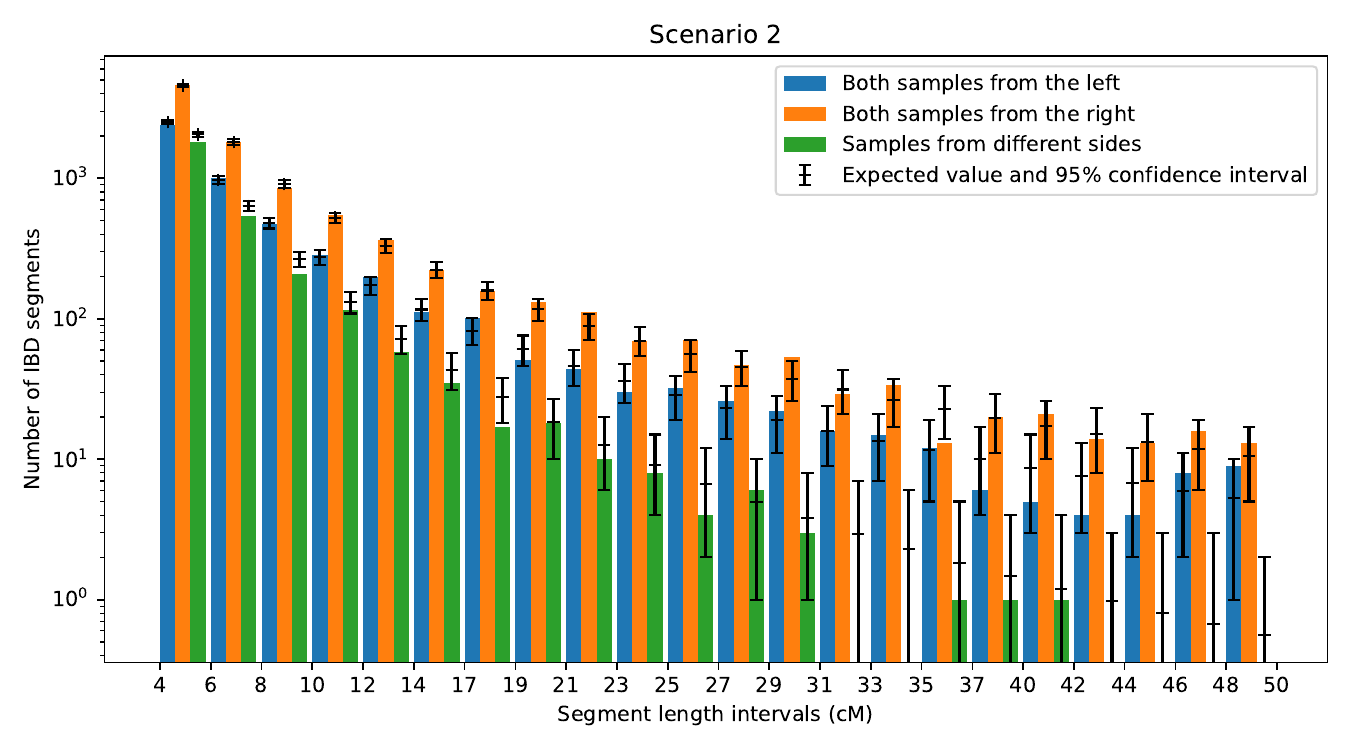}
    \caption{Number of IBD segments of different lengths observed in one simulation of the SLFV for scenarios 1 and 2 (see Table~\ref{table:simu_params}).
    The bar plot shows the sum of all observed IBD segments of each length interval for all pairs of individuals respectively both on the left side of the interface, both on the right side, and on different sides. The expected value computed with \eqref{eq:WM_E_N_geq_L} with the simulation parameters and the associated 95\% confidence interval of the Poisson distribution are shown in black for each length interval and each sampling configuration.}
    \label{supp:fig:observed_ibd}
\end{figure}

\begin{figure}[p]\ContinuedFloat
    \centering
    \includegraphics[width=0.9\textwidth]{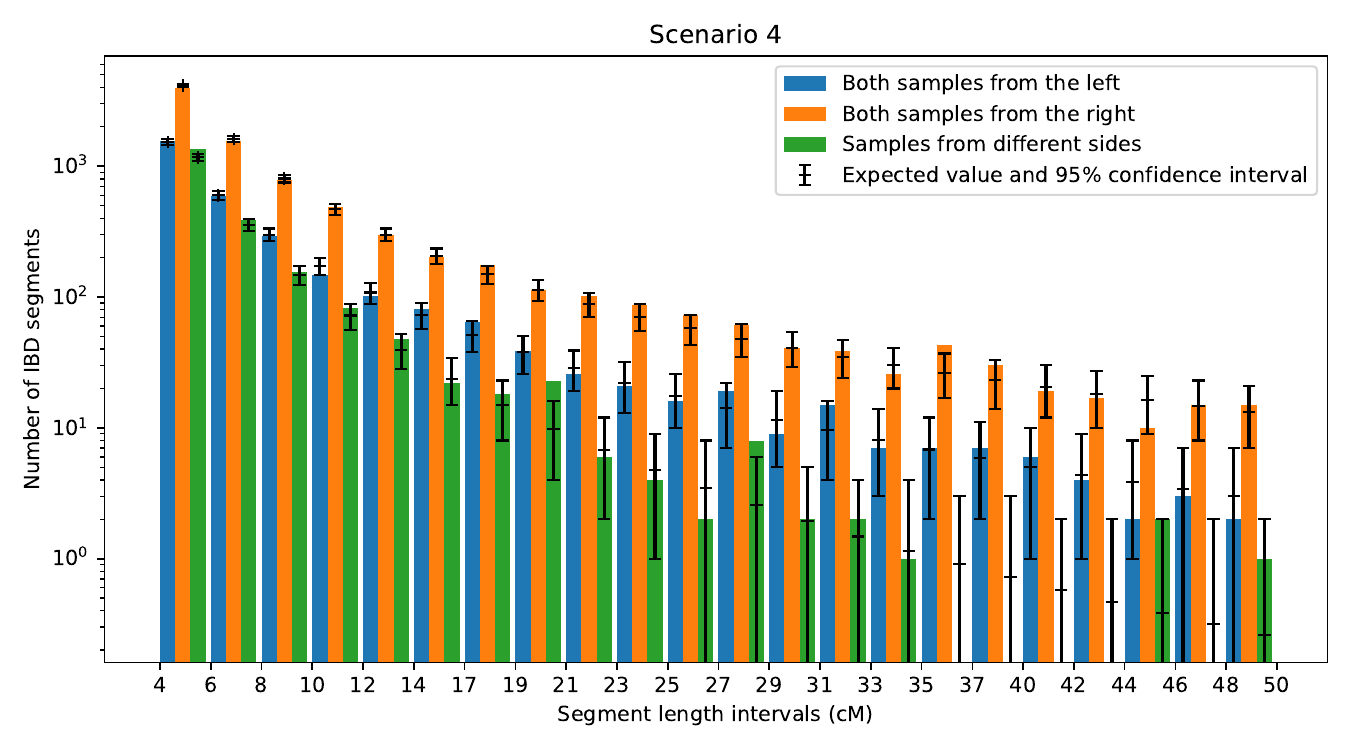}
    \caption{Corresponding graphic for scenario~4 (see Table~\ref{table:simu_params}).}
\end{figure}

\clearpage

\section{Proof of Lemma~\ref{lemma:stationary_dist_Y}} \label{sec:proof_lemma_h}

\begin{proof}[Proof of Lemma~\ref{lemma:stationary_dist_Y}]
    We start by proving that $(Y(t), t \geq 0)$ is ergodic with respect to the uniform probability measure on $[-R,R]$.
	Let $ \mathcal{L}^Y $ denote the infinitesimal generator of $ (Y(t), t \geq 0) $.
	For $ f$ and $ g $ two bounded and measurable functions on a set $ B \subset \R $, let
	\begin{equation*}
	\langle f, g \rangle_B := \int_{B} f(x) g(x) dx.
	\end{equation*}
	We want to show
	\begin{equation} \label{symmetry_LY}
	\langle \mathcal{L}^Y f, g \rangle_{[-R,R]} = \langle f, \mathcal{L}^Y g \rangle_{[-R, R]}.
	\end{equation}
	Recall the definition of the operator $ E $ in \eqref{definition_operator_E} and note that, $ \mathcal{L}^Y f = \mathcal{L} E f $ on $[-R,R]$, where $ \mathcal{L} $ is defined in \eqref{def:generator_single_lineage}.
	In addition, since $ \Phi(x,y) = \Phi(y,x) $, for any $ f, g \in L^2(\R) $,
	\begin{equation*}
	\langle \mathcal{L}f, g \rangle_\R = \langle f, \mathcal{L}g\rangle_\R.
	\end{equation*}
	However, $ Ef \notin L^2(\R) $.
	To circumvent this, for $ A \geq R $, define
	\begin{equation*}
		\Phi^A(x,y) := \left\lbrace
		\begin{aligned}
		\Phi(x,y) \hphantom{0} & \text{ if } \abs{x} \leq A \text{ and } \abs{y} \leq A, \\
		0 \hphantom{\Phi(x,y)} & \text{ otherwise.}
		\end{aligned}
		\right.
	\end{equation*}
	Further let $ (\xi^A_t, t \geq 0) $ be a random walk on $ \R $ with generator
	\begin{align} \label{def_LA}
		\mathcal{L}^A f(x) := \int_\R \Phi^A(x,y) (f(y) - f(x)) dy
	\end{align}
	which coincides with $ (\xi_t, t \geq 0) $ up to the random time $ T_A := \inf \{ t \geq 0 : \abs{\xi_t} > A \} $.
	Finally, for $ \abs{x} \leq A $, define
	\begin{equation} \label{def_EA}
		E^A f(x) := \E[x]{f(\xi^A_{T^A_0})}, \quad \text{ with } \quad T^A_0 := \inf \{ t \geq 0 : \abs{\xi^A_t} \leq r_+ \}.
	\end{equation}
	Then the operator $ \mathcal{L}^A $ is self-adjoint in $ L^2([-A,A]) $ and
	\begin{align}
		\mathcal{L}^A E^A f(x) = 0 & \qquad \text{ for } R < \abs{x} \leq A, \label{EA_1} \\
		E^A f(x) = f(x) & \qquad \text{ for } \abs{x} \leq R. \label{EA_2}
	\end{align}
	As a result, for $ f, g : [-R,R] \to \R $ bounded and measurable, 
	\begin{align*}
		\langle \mathcal{L}^A E^A f, E^A g \rangle_{[-A,A]} = \langle E^A f,\mathcal{L}^A E^A g\rangle_{[-A,A]}.
	\end{align*}
	By \eqref{EA_1} and \eqref{EA_2}, this is
	\begin{align*}
		\langle \mathcal{L}^A E^A f, g \rangle_{[-R,R]} = \langle f, \mathcal{L}^A E^A g\rangle_{[-R,R]}.
	\end{align*}
	It thus remains to let $ A \to \infty $.
	First note that, for $ A $ large enough, $ \mathcal{L}^A f(x) = \mathcal{L} f(x) $ for all $ x \in [-R,R] $.
	Furthermore, since $ T_A \to \infty $ as $ A \to \infty $ almost surely, $ \xi^A_{T^A_0} \to \xi_{\alpha(0)} $ as $ A \to \infty $ almost surely.
	Hence, by dominated convergence, for all $ x \in \R $ and for all bounded and measurable $ f $,
	\begin{align} \label{convergence_EA}
		E^A f(x) \cvgas{A} E f(x).
	\end{align}
	We thus obtain
	\begin{align*}
		\langle \mathcal{L} E f, g \rangle_{[-R,R]} = \langle f, \mathcal{L} E g\rangle_{[-R,R]},
	\end{align*}
	which is \eqref{symmetry_LY}.
	As a result the uniform measure on $ [-R,R] $ is invariant for $ Y $.
	The fact that $ (Y(t), t \geq 0) $ is ergodic then follows from the form of its generator, noting that $ \mathcal{L}^Y f \equiv 0 $ implies that $ f $ is almost everywhere constant.
\end{proof}

Before proving the second part of Lemma~\ref{lemma:stationary_dist_Y}, let us show the following.

\begin{lemma} \label{lemma:cancellation}
	For any $ f : [-R,R] \to \R $ bounded and measurable,
	\begin{align*}
		\int_{-\infty}^{R} \int_{R}^{+\infty} \Phi(x,y) \left( E f(x) - E f(y) \right) dy dx = 0
	\end{align*}
	and
	\begin{align*}
		\int_{-R}^{+\infty} \int_{-\infty}^{-R} \Phi(x,y) \left( E f(x) - E f(y) \right) dy dx = 0
	\end{align*}
\end{lemma}

\begin{proof}
	Recall the definition of $ \mathcal{L}^A $ and $ E^A $ in \eqref{def_LA} and \eqref{def_EA}.
	By \eqref{EA_1}, for any $ f : [-R,R] \to \R $ bounded and measurable,
	\begin{align*}
		\int_{R}^{A} \mathcal{L}^A E^A f(x) dx = 0.
	\end{align*}
	Since $ \mathcal{L}^A $ is self-adjoint on $ [-A,A] $,
	\begin{align*}
		0 = \int_{-A}^{A} \int_{-A}^{A} E^A f(x) \Phi^A(x,y) \left( \1{y > R} - \1{x > R} \right) dy dx.
	\end{align*}
	But
	\begin{align*}
		\1{y > R} - \1{x > R} = \1{y > R, x \leq R} - \1{y \leq R, x > R}.
	\end{align*}
	The above term is thus
	\begin{align*}
		\int_{-A}^{R} \int_{R}^{A} \Phi^A(x,y) E^A f(x) dy dx - \int_{R}^{A} \int_{-A}^{R} \Phi^A(x,y) E^A f(x) dy dx.
	\end{align*}
	Since $ \Phi^A(x,y) = \Phi^A(y,x) $, we obtain
	\begin{align*}
		\int_{-A}^{R} \int_{R}^{A} \Phi^A(x,y) \left( E^A f(x) - E^A f(y) \right) dy dx = 0.
	\end{align*}
	Letting $ A \to \infty $ and using \eqref{convergence_EA}, we obtain the first statement of Lemma~\ref{lemma:cancellation}.
	The second statement follows by a similar argument.
\end{proof}
	
We now finish the proof of Lemma~\ref{lemma:stationary_dist_Y}.

\begin{proof}[Conclusion of proof of Lemma~\ref{lemma:stationary_dist_Y}]
	Let us start by computing $ \frac{1}{2R} \int_{-R}^R h^+ (x) dx $.
	By the first part of Lemma~\ref{lemma:stationary_dist_Y} and \eqref{function_h},
	\begin{align*}
		\frac{1}{2R} \int_{-R}^R h^+ (x) dx = \frac{u}{2R} \int_{-R}^{R} \int_{-\infty}^{R} \Phi(x,y) \left( E \iota(y) - x \right) dy dx + \frac{u}{2R} \int_{-R}^{R} \int_{R}^{+\infty} \Phi(x,y) (y-x) dy dx.
	\end{align*}
	Since $ E\iota(y) = y $ when $ \abs{y} \leq R $ and $ \Phi(x,y) = 0 $ when $ \abs{x - y} > 2R $, this is
	\begin{multline*}
		\frac{1}{2R} \int_{-R}^R h^+ (x) dx = \frac{u}{2R} \int_{-R}^{R} \int_{-R}^{R} \Phi(x,y) \left( y-x \right) dy dx  + \frac{u}{2R} \int_{-R}^{\infty} \int_{-\infty}^{-R} \Phi(x,y) \left( E \iota(y) - E \iota (x) \right) dy dx \\ + \frac{u}{2R} \int_{-\infty}^{R} \int_{R}^{+\infty} \Phi(x,y) (y-x) dy dx.
	\end{multline*}
	The first term on the right hand side is zero because $ \Phi(x,y) = \Phi(y,x) $ and the second term is zero by Lemma~\ref{lemma:cancellation}.
	Replacing $ y $ by $ x + z $ in the last term, we have
	\begin{align*}
		\frac{1}{2R} \int_{-R}^R h^+ (x) dx = \frac{u}{2R} \int_{-\infty}^{R} \int_{0}^{+\infty} \Phi(x,x+z) \1{x + z > R} z dz dx.
	\end{align*}
	But for $ x+z > R $, $ B(x+z,r_-) \cap \mathbb{H}^- = \emptyset $ and
	\begin{equation*}
	\Phi(x, x+z) = \frac{\abs{B(x, R) \cap B(x+z,R)}}{V_{R}^2} = \Phi(R, R + z).
	\end{equation*}
	Furthermore, the expression above is zero when $ z \geq 2R $.
	Changing the order of integration, we obtain
	\begin{align*}
		\frac{1}{2R} \int_{-R}^R h^+ (x) dx &= \frac{u}{2R} \int_{0}^{2R} \Phi(R,R + z) z \int_{R - z}^{R} dx dz \\
		&= \frac{u}{2R} \int_{0}^{2R} \Phi(R,R + z) z^2 dz \\
		&= \frac{\sigma_+^2}{4 R}.
	\end{align*}
	By the same argument, one arrives at
	\begin{align*}
		\frac{1}{2R} \int_{-R}^R h^- (x) dx = \frac{\sigma_-^2}{4 R}.
	\end{align*}
\end{proof}

\bibliography{TPB_biblio}

\newcommand{\etalchar}[1]{$^{#1}$}
\begin{thebibliography}{CWWP14}

\bibitem[AFB{\etalchar{+}}17]{aguillon_deconstructing_2017}
Stepfanie~M. Aguillon, John~W. Fitzpatrick, Reed Bowman, Stephan~J. Schoech,
  Andrew~G. Clark, Graham Coop, and Nancy Chen.
\newblock Deconstructing isolation-by-distance: {{The}} genomic consequences of
  limited dispersal.
\newblock {\em PLOS Genetics}, 13(8):e1006911, 2017.

\bibitem[Ald78]{aldous_stopping_1978}
David Aldous.
\newblock Stopping times and tightness.
\newblock {\em The Annals of Probability}, 6(2):335--340, 1978.

\bibitem[APSN19]{al-asadi_estimating_2019}
Hussein {Al-Asadi}, Desislava Petkova, Matthew Stephens, and John Novembre.
\newblock Estimating recent migration and population size surfaces.
\newblock {\em PLOS Genet}, 15(1):365536, 2019.

\bibitem[BB11]{browning_fast_2011}
Brian~L. Browning and Sharon~R. Browning.
\newblock A fast, powerful method for detecting identity by descent.
\newblock {\em The American Journal of Human Genetics}, 88(2):173--182, 2011.

\bibitem[BB15]{browning_accurate_2015}
Sharon~R. Browning and Brian~L. Browning.
\newblock Accurate {{Non-parametric Estimation}} of {{Recent Effective
  Population Size}} from {{Segments}} of {{Identity}} by {{Descent}}.
\newblock {\em The American Journal of Human Genetics}, 97(3):404--418, 2015.

\bibitem[BB20]{browning_probabilistic_2020}
Sharon~R. Browning and Brian~L. Browning.
\newblock Probabilistic {{Estimation}} of {{Identity}} by {{Descent Segment
  Endpoints}} and {{Detection}} of {{Recent Selection}}.
\newblock {\em The American Journal of Human Genetics}, 107(5):895--910, 2020.

\bibitem[BBG{\etalchar{+}}16]{baharian_great_2016}
Soheil Baharian, Maxime Barakatt, Christopher~R. Gignoux, Suyash Shringarpure,
  Jacob Errington, William~J. Blot, Carlos~D. Bustamante, Eimear~E. Kenny,
  Scott~M. Williams, Melinda~C. Aldrich, et~al.
\newblock The great migration and {{African-American}} genomic diversity.
\newblock {\em PLoS genetics}, 12(5):e1006059, 2016.

\bibitem[BDE02]{barton_neutral_2002}
Nick~H. Barton, Frantz Depaulis, and Alison~M. Etheridge.
\newblock Neutral evolution in spatially continuous populations.
\newblock {\em Theoretical Population Biology}, 61(1):31--48, 2002.

\bibitem[BEKV13]{barton_inference_2013}
Nick~H. Barton, Alison~M. Etheridge, Jerome Kelleher, and Amandine V{\'e}ber.
\newblock Inference in two dimensions: Allele frequencies versus lengths of
  shared sequence blocks.
\newblock {\em Theoretical population biology}, 87:105--119, 2013.

\bibitem[BEV10]{barton_new_2010}
Nick~H. Barton, Alison~M. Etheridge, and Amandine V{\'e}ber.
\newblock A new model for evolution in a spatial continuum.
\newblock {\em Electronic Journal of Probability}, 15(7):162--216, 2010.

\bibitem[BGT89]{bingham_regular_1989}
Nicholas~H. Bingham, Charles~M. Goldie, and Jef~L. Teugels.
\newblock {\em Regular Variation}, volume~27 of {\em Encyclopedia of
  {{Mathematics}} and Its {{Applications}}}.
\newblock {Cambridge university press, Cambridge}, 1989.

\bibitem[BR19]{bradburd_spatial_2019}
Gideon~S. Bradburd and Peter~L. Ralph.
\newblock Spatial {{Population Genetics}}: {{It}}'s {{About Time}}.
\newblock {\em Annual Review of Ecology, Evolution, and Systematics},
  50(1):427--449, 2019.

\bibitem[BRK20]{battey_space_2020}
C.~J. Battey, Peter~L. Ralph, and Andrew~D. Kern.
\newblock Space is the {{Place}}: {{Effects}} of {{Continuous Spatial
  Structure}} on {{Analysis}} of {{Population Genetic Data}}.
\newblock {\em Genetics}, 215(1):193--214, 2020.

\bibitem[CHGG16]{coffman_computationally_2016}
Alec~J. Coffman, Ping~Hsun Hsieh, Simon Gravel, and Ryan~N. Gutenkunst.
\newblock Computationally {{Efficient Composite Likelihood Statistics}} for
  {{Demographic Inference}}.
\newblock {\em Molecular Biology and Evolution}, 33(2):591--593, 2016.

\bibitem[CRN16]{chiang_conflation_2016}
Charleston W~K Chiang, Peter Ralph, and John Novembre.
\newblock Conflation of {{Short Identity-by-Descent Segments Bias Their
  Inferred Length Distribution}}.
\newblock {\em G3 Genes|Genomes|Genetics}, 6(5):1287--1296, 2016.

\bibitem[CWWP14]{carmi_renewal_2014}
Shai Carmi, Peter~R. Wilton, John Wakeley, and Itsik Pe'er.
\newblock A renewal theory approach to {{IBD}} sharing.
\newblock {\em Theoretical Population Biology}, 97:35--48, 2014.

\bibitem[FRP22]{fournier_haplotype-based_2022}
Romain Fournier, David Reich, and Pier~Francesco Palamara.
\newblock Haplotype-based inference of recent effective population size in
  modern and ancient {{DNA}} samples.
\newblock {\em bioRxiv}, page 2022.08.03.501074, 2022.

\bibitem[GGW16]{guindon_demographic_2016}
Stephane Guindon, Hongbin Guo, and David Welch.
\newblock Demographic inference under the coalescent in a spatial continuum.
\newblock {\em Theoretical population biology}, 111:43--50, 2016.

\bibitem[GM97]{griffiths_ancestral_1997}
Robert~C. Griffiths and Paul Marjoram.
\newblock An ancestral recombination graph.
\newblock {\em Institute for Mathematics and its Applications}, 87:257, 1997.

\bibitem[HS81]{harrison_skew_1981}
John~Michael Harrison and Lawrence~A. Shepp.
\newblock On skew {{Brownian}} motion.
\newblock {\em The Annals of Probability}, 9(2):309--313, 1981.

\bibitem[Hud83]{hudson_properties_1983}
Richard~R. Hudson.
\newblock Properties of a neutral allele model with intragenic recombination.
\newblock {\em Theoretical population biology}, 23(2):183--201, 1983.

\bibitem[IM63]{ito_brownian_1963}
K.~It{\^o} and H.~P. McKean.
\newblock Brownian motions on a half line.
\newblock {\em Illinois Journal of Mathematics}, 7:181--231, 1963.

\bibitem[IP16]{iksanov_functional_2016}
Alexander Iksanov and Andrey Pilipenko.
\newblock A functional limit theorem for locally perturbed random walks.
\newblock {\em Probability and Mathematical Statistics}, 36(2):353--368, 2016.

\bibitem[KW64]{kimura_stepping_1964}
Motoo Kimura and George~H. Weiss.
\newblock The stepping stone model of population structure and the decrease of
  genetic correlation with distance.
\newblock {\em Genetics}, 49(4):561, 1964.

\bibitem[Lej06]{lejay_constructions_2006}
Antoine Lejay.
\newblock On the constructions of the skew {{Brownian}} motion.
\newblock {\em Probability Surveys}, 3:413--466, 2006.

\bibitem[L{\'e}p78]{lepingle_sur_1978}
Dominique L{\'e}pingle.
\newblock Sur le comportement asymptotique des martingales locales.
\newblock In {\em S\'eminaire de {{Probabilit\'es XII}}}, pages 148--161.
  {Springer}, 1978.

\bibitem[Lin88]{lindsay_composite_1988}
B.~G. Lindsay.
\newblock Composite likelihood methods.
\newblock {\em Comtemporary Mathematics}, 80(1):221--239, 1988.

\bibitem[LL10]{lawler_random_2010}
Gregory~F. Lawler and Vlada Limic.
\newblock {\em Random Walk: A Modern Introduction}, volume 123 of {\em
  Cambridge {{Studies}} in {{Advanced Mathematics}}}.
\newblock {Cambridge University Press, Cambridge}, 2010.

\bibitem[Mal48]{malecot_les_1948}
Gustave Mal{\'e}cot.
\newblock {\em Les {{Math\'ematiques}} de l'{{H\'er\'edit\'e}}}.
\newblock {Masson et Cie., Paris}, 1948.

\bibitem[Nag76]{nagylaki_clines_1976}
Thomas Nagylaki.
\newblock Clines with {{Variable Migration}}.
\newblock {\em Genetics}, 83(4):867--886, 1976.

\bibitem[NBK{\etalchar{+}}08]{nelson_population_2008}
Matthew~R. Nelson, Katarzyna Bryc, Karen~S. King, Amit Indap, Adam~R. Boyko,
  John Novembre, Linda~P. Briley, Yuka Maruyama, Dawn~M. Waterworth, G{\'e}rard
  Waeber, et~al.
\newblock The {{Population Reference Sample}}, {{POPRES}}: A resource for
  population, disease, and pharmacological genetics research.
\newblock {\em The American Journal of Human Genetics}, 83(3):347--358, 2008.

\bibitem[NGY{\etalchar{+}}16]{ni_probabilistic_2016}
Xumin Ni, Wei Guo, Kai Yuan, Xiong Yang, Zhiming Ma, Shuhua Xu, and Shihua
  Zhang.
\newblock A {{Probabilistic Method}} for {{Estimating}} the {{Sharing}} of
  {{Identity}} by {{Descent}} for {{Populations}} with {{Migration}}.
\newblock {\em IEEE/ACM Transactions on Computational Biology and
  Bioinformatics}, 13(2):281--290, 2016.

\bibitem[NJB{\etalchar{+}}08]{novembre_genes_2008}
John Novembre, Toby Johnson, Katarzyna Bryc, Zolt{\'a}n Kutalik, Adam~R. Boyko,
  Adam Auton, Amit Indap, Karen~S. King, Sven Bergmann, Matthew~R. Nelson,
  et~al.
\newblock Genes mirror geography within {{Europe}}.
\newblock {\em Nature}, 456(7218):98--101, 2008.

\bibitem[PLDP12]{palamara_length_2012}
Pier~Francesco Palamara, Todd Lencz, Ariel Darvasi, and Itsik Pe'er.
\newblock Length distributions of identity by descent reveal fine-scale
  demographic history.
\newblock {\em The American Journal of Human Genetics}, 91(5):809--822, 2012.

\bibitem[Por79a]{portenko_diffusion_1979}
N.~I. Portenko.
\newblock Diffusion processes with generalized drift coefficients.
\newblock {\em Theory of Probability \& Its Applications}, 24(1):62--78, 1979.

\bibitem[Por79b]{portenko_stochastic_1979}
N.~I. Portenko.
\newblock Stochastic differential equations with generalized drift vector.
\newblock {\em Theory of Probability \& Its Applications}, 24(2):332--347,
  1979.

\bibitem[PP13]{palamara_inference_2013}
Pier~Francesco Palamara and Itsik Pe'er.
\newblock Inference of historical migration rates via haplotype sharing.
\newblock {\em Bioinformatics}, 29(13):i180--i188, 2013.

\bibitem[RC13]{ralph_geography_2013}
Peter Ralph and Graham Coop.
\newblock The geography of recent genetic ancestry across {{Europe}}.
\newblock {\em PLoS Biol}, 11(5):e1001555, 2013.

\bibitem[RCB17]{ringbauer_inferring_2017}
Harald Ringbauer, Graham Coop, and Nick~H. Barton.
\newblock Inferring {{Recent Demography}} from {{Isolation}} by {{Distance}} of
  {{Long Shared Sequence Blocks}}.
\newblock {\em Genetics}, 205(3):1335--1351, 2017.

\bibitem[Rou97]{rousset_genetic_1997}
Fran{\c c}ois Rousset.
\newblock Genetic differentiation and estimation of gene flow from
  {{F-statistics}} under isolation by distance.
\newblock {\em Genetics}, 145(4):1219--1228, 1997.

\bibitem[RY13]{revuz_continuous_2013}
Daniel Revuz and Marc Yor.
\newblock {\em Continuous Martingales and {{Brownian}} Motion}, volume 293.
\newblock {Springer Science \& Business Media}, 2013.

\bibitem[Saw76]{sawyer_results_1976}
Stanley Sawyer.
\newblock Results for the {{Stepping Stone Model}} for {{Migration}} in
  {{Population Genetics}}.
\newblock {\em The Annals of Probability}, 4(5):699--728, 1976.

\bibitem[SHM00]{sharbel_genetic_2000}
Timothy~F. Sharbel, Bernhard Haubold, and Thomas {Mitchell-Olds}.
\newblock Genetic isolation by distance in {{Arabidopsis}} thaliana:
  Biogeography and postglacial colonization of {{Europe}}.
\newblock {\em Molecular Ecology}, 9(12):2109--2118, 2000.

\bibitem[Wal78]{walsh_diffusion_1978}
John~B. Walsh.
\newblock A diffusion with a discontinuous local time.
\newblock {\em Ast\'erisque}, 52(53):37--45, 1978.

\bibitem[Wri43]{wright_isolation_1943}
Sewall Wright.
\newblock Isolation by distance.
\newblock {\em Genetics}, 28(2):114, 1943.

\bibitem[ZBB20]{zhou_fast_2020}
Ying Zhou, Sharon~R. Browning, and Brian~L. Browning.
\newblock A {{Fast}} and {{Simple Method}} for {{Detecting Identity-by-Descent
  Segments}} in {{Large-Scale Data}}.
\newblock {\em The American Journal of Human Genetics}, 106(4):426--437, 2020.

\end{thebibliography}

\end{document}